\documentclass[12pt,letterpaper,twoside]{article}

\usepackage{amsmath,amsthm,amsfonts,amssymb}
\usepackage{times}
\usepackage{inputenc}
\usepackage{amsmath}
\usepackage{nicefrac}
\usepackage{amsthm}
\usepackage{amsfonts}
\usepackage{verbatim}
\usepackage{color}
\usepackage{bbm}
\usepackage{fancyhdr}
\usepackage[arrow, matrix, curve]{xy}

\newcommand{\pdftitle}{Path-by-path uniqueness of infinite-dimensional stochastic differential equations}

\newcommand{\pdfauthor}{Lukas Wresch}

\author{\pdfauthor\\\small{Faculty of Mathematics, Bielefeld University, Germany, E-mail: \href{mailto:wresch@math.uni-bielefeld.de}{wresch@math.uni-bielefeld.de}}}

\title{\pdftitle}

\usepackage[plainpages=false,pdfpagelabels=true,bookmarks=true,pdfauthor={\pdfauthor},pdftitle={\pdftitle}]{hyperref}

\newtheoremstyle{Theoremstyle} 
                        {1.5em}    
                        {2.5em}    
                        {}         
                        {}         
                        {\bfseries}
                        {}         
                        {\newline} 
                        {\raisebox{0.6em}{\thmname{#1}\thmnumber{ #2}\thmnote{ (#3)}}}

\newcounter{satz}[section]
\renewcommand{\thesatz}{\thesection.\arabic{satz}}
\newcounter{ref}[satz]

\makeatletter
\def\HyPsd@CatcodeWarning#1{}

\newcommand{\customlabel}[1]{%
	 \stepcounter{ref}%
   \protected@write \@auxout{}{\string\newlabel{#1}{{\thesatz.\arabic{ref}}{\thepage}{\thesatz.\arabic{ref}}{#1}{}}}%
	\hypertarget{#1}{}%
}

\newcommand{\xlabel}[1]{%
	 \label{#1}%
}
\makeatother

\theoremstyle{Theoremstyle}
\newtheorem{lem}[satz]{Lemma}
\newtheorem{pro}[satz]{Proposition}
\newtheorem{cor}[satz]{Corollary}
\newtheorem{re}[satz]{Remark}
\newtheorem{de}[satz]{Definition}

\newtheorem{thm}[satz]{Theorem}
\newtheorem{ass}[satz]{Assumption}
\newtheorem*{prof}{Proof}

\newcommand{\R}{\mathbb{R}}
\newcommand{\N}{\mathbb{N}}

\newcommand{\Z}{\mathbb{Z}}

\newcommand{\D}{\mathbb{D}}
\renewcommand{\P}{\mathbb{P}}

\newcommand{\eps}{\varepsilon}
\newcommand{\E}{\mathbb E}

\renewcommand{\d}{\,\mathrm{d}}

\newcommand{\bbm}{\mathbbm}
\newcommand{\ca}{\mathcal}
\newcommand{\mrm}{\mathrm}

\newcommand{\op}[1]{\operatorname{#1}}

\newcommand{\la}{\langle}
\newcommand{\ra}{\rangle}

\newcommand{\h}{\hat}

\newcommand{\bfrac}[2]{\genfrac{}{}{0pt}{}{#1}{#2}}

\begin{document}

\maketitle

\begin{abstract}
\noindent
Consider the stochastic differential equation $\mathrm dX_t = -A X_t \,\mathrm dt + f(t, X_t) \,\mathrm dt  +  \mathrm dB_t$ in a (possibly infinite-dimensional) separable Hilbert space, where $B$ is a cylindrical Brownian motion and $f$ is a just measurable, bounded function. If the components of $f$ decay to 0 in a faster than exponential way we establish path-by-path uniqueness for mild solutions of this stochastic differential equation. This extends A.~M.~Davie's result from $\R^d$ to Hilbert space-valued stochastic differential equations.
\end{abstract}

\section{Preliminaries}

\subsection{Framework \& Main result}

Let us consider the following stochastic differential equation (SDE)

\[ \begin{cases} \d \!\!\!\!\!\!& x(t) = -A x(t) \,\mrm dt + f(t, x(t)) \,\mrm dt + \mrm dB_t \\ \!\!\!\!\!\!& x(0)\! = x_0 . \end{cases}  \tag{SDE} \]

on a separable Hilbert space $H$ in mild form i.e.~a solution $x$ satisfies

\[ x(t) = e^{-tA} x_0 + \int\limits_0^t e^{-(t-s)A} f(s, x(s)) \d s  +  \int\limits_0^t  e^{- (t-s)A} \d B_s  \qquad \P\text{-a.s.},\ \forall t\in [0,T] .  \tag{1.1.1}  \]

Let, as in the previous article \cite{Wre16}, $H$ be a separable Hilbert space over $\R$.
Let\\$(\Omega, \ca F, (\ca F_t)_{t\in[0,\infty[}, \P, (B_t)_{t\in[0,\infty[} )$ be a filtered stochastic basis with sigma-algebra $\ca F$, a right-continuous, normal filtration $\ca F_t \subseteq \ca F$, a probability measure $\P$ and $(B_t)_{t\in[0,\infty[}$ an $\ca F_t$-Brownian motion on $(\Omega, \ca F, \P)$ taking values in $\R^\N$.
Let $A \colon D(A) \longrightarrow H$ be a positive definite, self-adjoint, linear operator with trivial kernel such that $A^{-1}$ is trace-class. Hence, there exists an orthonormal basis $(e_n)_{n\in\N}$ of $H$ such that

\[ \hspace{28mm}  A e_n = \lambda_n e_n, \qquad\qquad  \lambda_n > 0, \ \forall n \in\N  \]

with

\[ \hspace{31.5mm}  \lambda_n \leq \lambda_{n+1} , \qquad\qquad  \hspace{14.8mm}  \forall n \in\N . \]

By fixing this basis $(e_n)_{n\in\N}$ we identify $H$ with $\ell^2$, so that $H \cong \ell^2 \subseteq \R^\N$.
Let $f \colon [0,1] \times H \longrightarrow H$ be a bounded, Borel measurable map.

\begin{re}[Existence of weak solutions]
\label{RE-EXISTENCE}

Using Girsanov's Theorem (see e.g.~\cite[Theorem I.0.2]{LR15}) we can construct a filtered stochastic basis as above and an $(\ca F_t)_{t\in[0,\infty[}$-adapted stochastic process $(X_t)_{t\in[0,T[}$ with $\P$-a.s.~continuous sample paths in $H$ which solves (SDE). I.e.~we have

\[ \begin{cases} \d \!\!\!\!\!\!& X_t = -A X_t \,\mrm dt + f(t, X_t) \,\mrm dt + \mrm dB_t \\ \!\!\!\!\!\!& X_0\! = x_0 . \end{cases}  \]

\end{re}

On an \textit{arbitrary} filtered stochastic basis $(\Omega, \ca F, (\ca F_t)_{t\in[0,\infty[}, \P, (B_t)_{t\in[0,\infty[} )$, as above, for which a priori it is not clear whether it carries a solution $(X_t)_{t\in[0,T[}$ as in Remark \ref{RE-EXISTENCE}, we study the equation (SDE)

\[ \begin{cases} \d \!\!\!\!\!\!& x(t) = -A x(t) \,\mrm dt + f(t, x(t)) \,\mrm dt + \mrm dB_t \\ \!\!\!\!\!\!& x(0)\! = x_0 . \end{cases}  \tag{SDE} \]

for bounded, measurable $f \colon [0,T] \times H \longrightarrow H$ with $T>0$ and $x_0\in H$.
We consider the so-called path-by-path approach
where equation (SDE) is not considered as a stochastic differential equation, but as a random integral equation in the mild sense. More precisely, in the path-by-path picture we first plug in an $\omega\in\Omega$ into the corresponding integral equation (IE) of the mild form of equation (SDE)

\customlabel{SDE}
\[ x_t = e^{-tA} x_0 + \int\limits_0^t e^{-(t-s)A} f(s, x_s) \d s + \left(\int\limits_0^t  e^{- (t-s)A} \d B_s\right) (\omega)  \tag*{$\text{(IE)}_\omega$} \]

and aim to find a (unique) continuous function $x \colon [0,T] \longrightarrow H$ satisfying this equation, which can now be considered as an ordinary integral equation (IE), that is perturbed by an Ornstein--Uhlenbeck path $Z^A(\omega)$. If such a (unique) function can be found for almost all $\omega\in\Omega$, the map $\omega \longmapsto x$ is called a (unique) path-by-path solution to the equation (SDE).
Naturally, this notion of uniqueness is much stronger than the usual pathwise uniqueness considered in the theory of SDEs.\\

The main result of this article states that on every filtered stochastic basis as above there exists a \textit{unique} mild solution to the equation (SDE) in the path-by-path sense. Although, in the finite dimensional setting many papers have been written about path-by-path uniqueness (see for example \cite{Dav07}, \cite{Sha14}, \cite{BFGM14}, \cite{Pri15}) to the best of our knowledge this is the first result in a general infinite-dimensional Hilbert space setting. However, for the special case where $H = L^2([0,1], \R)$ and $A=\Delta$ path-by-path uniqueness has been shown recently in \cite{BM16} for space-time white noise.\\

Let us now state the assumptions on the drift $f$ and the main result.

\begin{ass}
\label{ASS}

From now on let $f \colon [0,1] \times H \longrightarrow H$ be a Borel measurable map with components $f = f^{(n)}$ w.r.t.~our fixed basis $(e_n)_{n\in\N}$ satisfying the following conditions

\begin{align*}
\|f\|_H &= \sup\limits_{t\in[0,1], x\in H} |f(t,x)|_H \leq 1 ,  \\
\|f\|_{\infty,A} :&= \sup\limits_{t\in[0,1], x\in H} \sum\limits_{n\in\N} \lambda_n e^{2\lambda_n} |f^{(n)}(t, x)|^2 \leq 1  
\intertext{and}
\|f^{(n)}\|_\infty &= \sup\limits_{t\in[0,1], x\in H} |f^{(n)}(t,x)| \leq \exp\left( - e^{n^\gamma} \right)
\end{align*}

for some $\gamma > 6$.

\end{ass}

\begin{thm}[Main result]
\xlabel{THM-MAIN}

Let $A$ and $f$ be as above and assume that $f$ fulfills Assumption \ref{ASS}.
Given \textit{any} filtered stochastic basis $(\Omega, \ca F, (\ca F_t)_{t\in[0,\infty[}, \P, (B_t)_{t\in[0,\infty[})$
there exists $\Omega_0\in \ca F$ with $\P[\Omega_0]=1$ such that for every $\omega\in\Omega_0$ we have

\[ \# \{ g \in \ca C( [0,T], H) | g \text{ solves } \text{(IE)}_\omega \} = 1 ,  \]

i.e.~(SDE) has a path-by-path unique mild solution.

\end{thm}

Theorem \ref{THM-MAIN} follows from the following

\begin{pro}
\xlabel{PRO-MAIN}

Let $A$ and $f$ be as in Theorem \ref{THM-MAIN}.
Let $(\Omega, \ca F, (\ca F_t)_{t\in[0,\infty[}, \P, (B_t)_{t\in[0,\infty[})$ be a filtered stochastic basis and $(X_t)_{t\in[0,\infty[}$ a solution of (SDE) (as in Remark \ref{RE-EXISTENCE}).
Then path-by-path uniqueness holds, i.e.~there exists $\Omega_0\in\ca F$ with $\P[\Omega_0]=1$ such that

\[ \# \{ g \in \ca C( [0,T], H) | g \text{ solves } \text{(IE)}_\omega \} = 1  \]

holds for every $\omega\in\Omega_0$.

\end{pro}

\begin{prof}[of Theorem \ref{THM-MAIN}]

Take an arbitrary filtered probability space and let $( (X^{1}_t)_{t\in[0,\infty[}, (B_t)_{t\in[0,\infty[} )$ and\\$( (X^{2}_t)_{t\in[0,\infty[}, (B_t)_{t\in[0,\infty[} )$ be two weak solutions driven by the same cylindrical $(\ca F_t)_{t\in[0,\infty[}$-Brownian motion.
Then by Proposition \ref{PRO-MAIN} it follows that path-by-path uniqueness, and hence pathwise uniqueness, holds i.e.~$X^1=X^2$ $\P$-a.s.
Hence the Yamada--Watanabe Theorem (see \cite[Theorem 2.1]{RSZ08}) implies that there exists even a \textit{strong} solution for equation (SDE).
In conclusion, by invoking Proposition \ref{PRO-MAIN} again, this proves the existence \textit{and} path-by-path uniqueness of solutions on \textit{every} filtered stochastic basis $(\Omega, \ca F, (\ca F_t)_{t\in[0,\infty[}, \P, (B_t)_{t\in[0,\infty[})$.

\end{prof}

\begin{re}

Set $\Omega := L^2( [0,T], H)$ and $\P$ such that the projection $\pi_t(\omega) := \omega(t)$ is a cylindrical Brownian motion.
As in the introduction consider the map

\[ Z^A\colon L^2( [0,T], H) \longrightarrow \ca C([0,T], H), \qquad \omega \longmapsto \left( t \mapsto \int\limits_0^t e^{-(t-s)A} \d\omega(s) \right)  . \]

Note that due to \cite[Theorem 5.2]{DZ92} $\left(\P \circ Z^A\right)^{-1}$ equals $N(0,K)$, the Gaussian measure on $L^2( [0,T], H)$ with covariance operator $K$ defined by

\[ (K \varphi)(t)  =  \int\limits_0^T k(t,s) \varphi(s) \d s , \]
where
\[ k(t,s) = \int\limits_0^{t \wedge s}  e^{-(t-r)A} \left(e^{-(s-r)A}\right)^\star \d r  \]

and $N(0,K)[ Z^A(\Omega) ] = 1$. Note that, since $Z^A$ is injective, Kuratowski's Theorem (see \cite[Theorem A1.7]{Kal97}) implies that $Z^A(\Omega)$ is a Borel set.

\medskip

Let $f$ be as in Assumption \ref{ASS} then path-by-path uniqueness holds for the SDE

\[ \mathrm d x_t = -Ax_t \mathrm dt + f(t, x_t) \mathrm dt + \mathrm \omega(t) .  \]

I.e.~there exists $\Omega_0\subseteq \ca C([0,T], H)$ with $\P[\Omega_0]=1$ such that for every $\omega\in\Omega_0$ there exists a unique function $g \in \ca C( [0,T], H)$ solving the above equation.

\end{re}

\begin{re}

If the function $f$ is independent of time and $A^{-1+\delta}$ is trace class for some $\delta \in ]0,1[$ is has already been proven that a strong mild solution to (SDE) exists for $\mu$-a.a.~initial condition $x_0\in H$, where $\mu$ is the invariant measure of the Ornstein--Uhlenbeck process $Z^A$ (see \cite{DFPR13}). If in addition $f$ satisfies Assumption \ref{ASS}, our results improves this result to all initial conditions $x_0\in H$.

\end{re}

\subsection{Structure of the article \& Roadmap for the proof}

The structure of this article is the following: In the following section 2 we introduce approximation lattices and the notion of the effective dimension of an (infinite-dimensional) set. This is reminiscent to the Kolmogorov $\eps$-entropy, which was used in the proof of A.~V.~Shaposhniko (see \cite{Sha14}) for the finite-dimensional case. In the third section we prove two regularization by noise estimates of the map

\[ \varphi_{n,k} \colon (x,y) \longmapsto \int\limits_{k2^{-n}}^{(k+1)2^{-n}}  \left[ f(s, Z^A_s + x) - f(s, Z^A_s + y) \right]  \,\mathrm ds  , \]

which are based on the estimates previously obtained by the author in \cite{Wre16}. We show that for every $\delta > 0$

\[ |\varphi_{n,k}(x,y)|_H \leq C_\delta \left( \sqrt n 2^{-n/6} |x-y|_H + \eps_n \right) \]

with $\eps_n \overset{n\rightarrow\infty}\longrightarrow 0$ for all $\omega\in\Omega$ outside a set of mass $\delta$. Here, $x$ and $y$ are in an approximation lattice of a suitable subset $Q$ of $H$ which includes the image of $f$. For fixed $\omega\in\Omega$ the map $\varphi_{n,k}$ is therefore ``close to'' being Lipschitz continuous. This estimate acts as a replacement for the lack of regularity for the non-linearity $f$ in equation (SDE).\\

In the fourth section we extend these estimates: For sequences of functions $(h_m)_{m\in\N}$ converging to $h$ we prove, despite the lack of continuity in $f$, that

\[  \int\limits_0^1 f(s, Z^A_s + h_m(s)) \d s \overset{m\rightarrow\infty}\longrightarrow \int\limits_0^1 f(s, Z^A_s + h(s)) \d s  . \qquad \P\text{-a.s.}  \tag{1.2.1}  \]

This approximation theorem (Theorem \ref{THM-APPROX}) implies that the above map $\varphi_{n,k}$ is continuous and therefore enables us to extend the estimates of the previous section from an approximation lattice to all $x,y$ of $Q$ (Corollary \ref{COR-SIGMA-RHO}). The result obtained in this section is also necessary to justify the limiting argument in the proof of Theorem \ref{THM-FINAL}.\\

It turns out that in the proof of the main result (Theorem \ref{THM-MAIN}) we have to consider terms of type

\[ \sum\limits_{q=1}^N |\varphi_{n,k+q}(x_{q+1}, x_q)|_H   \tag{1.2.2}  \]

for a sequence of points $\{x_q \in Q | q=1, ..., N\}$. Using just the estimates of Section 3 for each term under the sum of (1.2.2) is, unfortunately, insufficient to prove the main result (Theorem \ref{THM-MAIN}) as this would merely give us an estimate of order $\ca O(2^{-(1/2-\eps)n} N)$.\\

To overcome this, in Section 5 we use the fact that the above points $x_q\in Q$ are values of a solution of an integral equation and hence can be well approximated by a one-step Euler approximation.
This enables us to prove much stronger estimates for expression (1.2.2). Namely, bounds of order $\ca O(2^{-n} N)$ (Theorem \ref{THM-EULER}).\\

Section 6 contains the proof of the main result (Theorem \ref{THM-MAIN}). As a first step of the proof of the main theorem we reduce the problem via Girsanov's Theorem to the following Proposition.

\begin{pro}[Reduction via Girsanov's Theorem]
\xlabel{PRO-GIRSANOV}

For every $f \colon [0,T] \times H \longrightarrow H$ be a Borel measurable function fulfilling Assumption \ref{ASS}.
Assume that for every process $(\tilde Z^A_t)_{t\in[0,\infty[}$ on $(\Omega, \ca F, (\ca F_t)_{t\in[0,\infty[})$ with $\tilde Z^A_0=0$, which is an Ornstein--Uhlenbeck process with drift term $A$ w.r.t.~some measure $\tilde\P$ on $(\Omega, \ca F)$, there exists a set $\Omega'_{\tilde Z^A} \subseteq \Omega$ with $\tilde\P[\Omega'_{\tilde Z^A}] = 1$ such that for all fixed $\omega\in \Omega'_{\tilde Z^A}$ the only function $u \in \ca C( [0,T], H)$ solving

\customlabel{DE-U}
\[ u(t) = \int\limits_0^t e^{-(t-s)A} \left( f(s, \tilde Z^A_s(\omega) + u(s)) - f(s, \tilde Z^A_s(\omega)) \right) \d s    \tag{\ref{DE-U}}  \]

for all $t\in[0,T]$ is the trivial solution $u \equiv 0$, then the assertion of Proposition \ref{PRO-MAIN} holds with $\Omega_0 := \Omega'_{\tilde Z^A}$, where $\tilde Z^A_t := X_t - e^{-tA} x_0$ with $X$ being a solution of (SDE). Recall that $X$ is an Ornstein--Uhlenbeck process under a measure $\tilde\P$ obtained via Girsanov transformation.

\end{pro}

\begin{re}[Dependence of $\Omega_0$]

The set of ``good omegas'' $\Omega_0$ of the main result \ref{THM-MAIN} therefore depends solely on the strong solution $X$, the initial condition $x_0$ and the drift $f$.

\end{re}

A proof of this proposition will be given in this section below.
Now, let $u$ be a function solving equation \eqref{DE-U} and let us write $\varphi_{n,k}(x) := \varphi_{n,k}(x,0)$.
To show that every solution to \eqref{DE-U} is trivial we use a discrete logarithmic Gronwall inequality of the form

\[ |u((k+1)2^{-n})|_H \leq |u(k2^{-n})|_H \left(1+C 2^{-n} \log( 1 / |u(k2^{-n})|_H ) \right) . \]

In Section 6 we first show that

\[ | u((k+1)2^{-n})  -  u(k2^{-n}) |_H  \approx  |\varphi_{n,k}(u(\cdot))|_H . \]

Subsequently, we construct functions $u_\ell \overset{\ell\rightarrow\infty}\longrightarrow u$, which are constant on the dyadic intervals $[k2^{-\ell}, (k+1)2^{-\ell}[$. Using the equation (1.2.1) mentioned above this can be rewritten as

\[ \lim\limits_{\ell\rightarrow\infty} |\varphi_{n,k}(u_\ell(\cdot))|_H   \leq   |\varphi_{n,k}(u_n(\cdot))|_H  +  \sum\limits_{\ell=n}^\infty |\varphi_{n,k}( u_{\ell+1}(\cdot))  ,  u_\ell(\cdot)  )|_H  .  \]

Splitting the integrals and using that $u_\ell$ is constant on dyadic intervals of size $2^{-\ell}$ we can bring this is in the somewhat more complicated form

\[   |\varphi_{n,k}(u(k2^{-n}))|_H  +  \sum\limits_{\ell=n}^\infty \sum\limits_{r=k2^{\ell+1-n}}^{(k+1)2^{\ell+1-n}} |\varphi_{\ell,r}( u((r+1)2^{-\ell-1})   , u(r2^{-\ell-1}) )|_H  .  \]

Using the estimates for $\varphi_{n,k}$ and expression (1.2.2) developed in the previous section we ultimately obtain an estimate of order

\[ | u((k+1)2^{-n})  -  u(k2^{-n}) |_H  \leq  C 2^{-n} |u(k2^{-n})|_H \log( 1 / |u(k2^{-n})|_H ) ,  \]

where we have to impose the somewhat technical condition that $0< |u(k2^{-n})|_H<1$. We therefore obtain a discrete $\log$-Type Gronwall inequality of the form

\[ |u((k+1)2^{-n})|_H \leq |u(k2^{-n})|_H \left(1+C 2^{-n} \log( 1 / |u(k2^{-n})|_H ) \right) , \]

which, similar to the standard Grownall Inequality, implies that $u$ has to be trivial (Corollary \ref{COR-FINAL}), so that the condition of Proposition \ref{PRO-GIRSANOV} is fulfilled completing the proof.

\begin{prof}[of Proposition \ref{PRO-GIRSANOV}]

Let $(X_t)_{t\in[0,T]}$ be a solution to (SDE). We set $\tilde Z^A_t := X_t - e^{-tA} x_0$ so that $\tilde Z^A$ is an Ornstein--Uhlenbeck process with drift term $A$ starting in $0$ under a measure $\tilde\P\approx\P$ obtained by Girsanov's Theorem as mentioned in Remark \ref{RE-EXISTENCE}.

Then, by assumption there is a set $\Omega'_{\tilde Z^A}$ with $\P[\Omega'_{\tilde Z^A}] = \tilde\P[\Omega'_{\tilde Z^A}] = 1$ such that for all $\omega \in \Omega'_{\tilde Z^A}$ every solution $u$ to equation \eqref{DE-U} is trivial.

\newpage

Let $\omega\in\Omega'_{\tilde Z^A}$ and $x \in \ca C( [0,T], H )$ be a solution to $\text{(IE)}_\omega$. We then have

\[ x_t = e^{-tA} x_0 + \int\limits_0^t e^{-(t-s)A} f(s, x_s) \d s  + \left(\int\limits_0^t  e^{- (t-s)A} \d B_s\right) (\omega)   .  \]

Setting $u_t := x_t - X_t(\omega)$ yields that

\begin{align*}
 u_t &= \int\limits_0^t e^{-(t-s)A} f(s, x_s) \d s - \int\limits_0^t e^{-(t-s)A} f(s, X_s(\omega)) \d s  \\
  &= \int\limits_0^t e^{-(t-s)A} ( f(s, u_s + X_s(\omega))  -   f(s, X_s(\omega)) ) \d s .
\end{align*}

By plugging in the definition of $\tilde Z^A$ and by setting

\[ \tilde f_{x_0}(t, z) := f(t, z + e^{-tA} x_0)   \]

we rewrite the above equation to

\[ u_t = \int\limits_0^t e^{-(t-s)A} ( \tilde f_{x_0}(s, u_s + \tilde Z^A_s(\omega))  -  \tilde f_{x_0}(s, \tilde Z^A_s(\omega) ) ) \d s \]

Since $\tilde Z^A$ is an Ornstein--Uhlenbeck process under $\tilde \P$ starting at zero and $\omega \in \Omega'_{\tilde Z^A}$ we conclude that $u \equiv 0$ and henceforth $x_t = X_t(\omega)$. Analogously, we obtain for any other solution $x'$ that $x'_t = X_t(\omega) = x_t$ so that all solutions of $\text{(IE)}_\omega$ coincide on $\Omega'_{\tilde Z^A}$ and are therefore unique.

\qed

\end{prof}

\section{Approximation Lattices}

In this section we define the set $Q$, where the function $u$ (see equation \eqref{DE-U} of Proposition \ref{PRO-GIRSANOV}) takes values in.
Additionally, we define the so-called \textit{effective dimension} of a set, which is a variant of the Kolmogorov $\eps$-entropy for lattices.
At the end of this section we estimate the effective dimension of our set $Q$.

\begin{de}[The set $Q$]
\xlabel{DE-Q}

We define

\[ Q := \{ x \in \R^\N  \colon |x|_\infty \leq 2, \ |x_n| \leq 2 \exp\left( - e^{n^\gamma} \right)  , \ x = (x_n)_{n \in \N} \}  , \]

where $\gamma$ is the constant from Assumption \ref{ASS}. Additionally, for $r \in\N$ we set

\[ Q_r := \{ x \in Q  \colon |x|_\infty \leq 2\cdot2^{-r} \} ,  \]

so that $Q_0=Q$.
Note that for $m\in\N$ the lattice $Q \cap 2^{-m} \Z^\N$ is the set of all points $x \in Q$, where the components $x_n$ of $x$ can be written as

\[  x_n = k_n 2^{-m} \]

with certain $k_n \in \Z$ for every $n\in\N$.

\end{de}

\begin{de}[Effective dimension]
\xlabel{DE-EFFDIM}

Let $B \subseteq \R^\N$ with $0\in B$. For points $x\in B$ we write $(x_n)_{n\in\N} = x$ for the components of $x$. For every $m\in\N$ we set

\[ d_m(B) := \sup\limits_{x \in B \cap 2^{-m} \Z^\N}  \inf\limits \left\{ \left. n \right|  x_{n'}=0 \ \forall n' \geq n \right\}   \in \bar\N := \N \cup \{ \infty \} . \]

I.e.~given any point $(x_n)_{n\in\N}$ in the set $B \cap 2^{-m} \Z^\N$, all components $x_n$ are zero for $n \geq d_m$ and $d_m$ is the smallest integer with this property.

\medskip

We define the \textit{effective dimension} of a set $B \subseteq \R^\N$ by
\[ \op{ed}\colon \{ B \subseteq \R^\N | 0 \in B \}  \longrightarrow {\bar \N}^\N  \]
\[ B  \longmapsto \op{ed}(B) := (d_m(B))_{m\in\N} .  \]

$B$ is called \textit{effectively finite-dimensional} if
\[ \op{ed}(B)_m < \infty , \qquad \forall m\in\N .  \]

Let $|\cdot|_1$ and $|\cdot|_2$ be two norm on $B$. $|\cdot|_1$ and $|\cdot|_2$ are called \textit{effectively equivalent} if for every $m\in\N$ they are equivalent on the restricted domain $B \cap 2^{-m} \Z^\N$. I.e.~for every $m\in\N$ there exists constants $c_m,C_m\in \R$ such that

\[ c_m|x|_1 \leq |x|_2 \leq C_m |x|_1  ,  \qquad  \forall x \in B \cap 2^{-m} \Z^\N  .  \]

\end{de}

\begin{pro}
\xlabel{EQUIVALENCE}

Let $B \subseteq \R^\N$ with $0\in B$ be an effectively finite-dimensional set then the norm $|\cdot|_2$ and the maximum norm $|\cdot|_\infty$ are effectively equivalent. More precisely, we have

\[ |x|_2 \leq \sqrt{\op{ed}(B)_m} |x|_\infty   , \qquad m \in \N, \ x \in B \cap 2^{-m} \Z^\N \]
and
\[ |x|_\infty \leq |x|_2                       , \qquad\qquad\qquad\quad\!\! m \in \N, \ x \in B \cap 2^{-m} \Z^\N . \]

\end{pro}

\begin{prof}

Let $m\in\N$. For every $x\in B \cap 2^{-m} \Z^\N$ we have

\[ |x|_2^2 = \sum\limits_{n=1}^\infty |x_n|^2 = \sum\limits_{n=1}^{\op{ed}(B)_m} |x_n|^2 \leq \op{ed}(H)_m |x|_\infty^2 \]

and

\[ |x|_\infty^2 \leq \sum\limits_{n=1}^\infty |x_n|^2 = |x|_2 .  \]

\qed

\end{prof}

\begin{lem}
\xlabel{EFFDIM}

For $r,m\in\N$ with $m\geq r$ we have

\[ \op{ed}(Q_r)_m \leq (\ln( m + 1 ))^{1/\gamma} .  \]

Note that this implies that $Q_r$ is effectively finite-dimensional for every $r\in\N$.

\end{lem}

\begin{prof}

Let $x \in Q_r \cap 2^{-m} \Z^\N$. Observe that every component $x_n$ is of the form $x_n = k_n 2^{-m}$ with
\[ k_n \in \{ - 2 \cdot 2^{m-r} , ..., 2 \cdot 2^{m-r} \} . \]

Set

\[ d_m := (\ln( m + 1 ))^{1/\gamma} .  \]

We are going to show that $k_n = 0$ holds for every $n \geq d_m$.

\[ |k_n| 2^{-m} = |x_n| \leq 2 \exp\left(- e^{n^\gamma} \right) \Rightarrow |k_n|  \leq 2^{m+1} \exp\left(- e^{n^\gamma} \right) ,  \]

which implies that

\[ |k_n| \leq 2^{m+1} \exp\left(- e^{n^\gamma} \right) \leq e^{\ln(2) (m+1)} \exp\left(- \exp\left( (d_m)^\gamma \right) \right) = e^{\ln(2) (m+1) - \exp\left( (d_m)^\gamma \right)} \]

\[ \hspace{33mm}  = e^{\ln(2) (m+1) - (m+1)} = e^{(\ln(2)-1) (m+1) } \leq e^{\ln(2)-1}  < 1 .  \]

In conclusion, $|k_n| = 0$ for all $n \geq d_m$ and hence we have

\[ \op{ed}(Q_r)_m \leq d_m = (\ln( m + 1 ))^{1/\gamma} .  \]

\qed

\end{prof}

\begin{thm}
\xlabel{KOLTIK}

Let $r\in\N$ and $m\in\N$. The number of points in the $m$-lattice of $Q_r$ can be estimated as follows

\[ \# (Q_r \cap 2^{-m} \Z^\N) \leq (4 \cdot 2^{m-r} + 1)^{\op{ed}(Q_r)_m}   \]

and

\[ \# (2Q_r \cap 2^{-m} \Z^\N) \leq (8 \cdot 2^{m-r} + 1)^{\op{ed}(2Q_r)_m} .   \]

\end{thm}

\begin{prof}

Let $m\in\N$ and $x \in Q_r \cap 2^{-m} \Z^\N$ and note that, as in the last proof, every component $x_n$ is of the form $x_n = k_n 2^{-m}$ with
\[ k_n \in \{ - 2 \cdot 2^{m-r} , ..., 2 \cdot 2^{m-r} \} . \]

$k_n$ can take at most $4 \cdot 2^{m-r} + 1$ different values in the dimensions $1 \leq n < \op{ed}(Q_r)_m$, so that the total number of points $x \in Q_r \cap 2^{-m} \Z^\N$ can be estimated by

\[ (4 \cdot 2^{m-r} + 1)^{\op{ed}(Q_r)_m}  . \]

Note that $k_n=0$ for $n \geq \op{ed}(Q_r)_m$. The second part of the assertion follows analogously.

\qed

\end{prof}

\begin{cor}
\xlabel{RE-PI}

Let $r\in\N$. For every $m\in\N$ there exists a map

\[ \pi_m^{(r)} \colon Q_r \longrightarrow Q_r \cap 2^{-m} \Z^\N \]

with the property that

\[ | x - \pi_m^{(r)}(x) |_\infty \leq 2^{-m}  \]

and

\[ \not{\exists} y \in Q_r \cap 2^{-m} \Z^\N \colon | x - y |_\infty < | x - \pi_m^{(r)}(x)|_\infty  \]

holds for all $x\in Q_r$, $m\in\N$ and $r\in\N$.

\end{cor}

\begin{prof}

Let $r, m\in\N$. By Theorem \ref{KOLTIK} and Lemma \ref{EFFDIM} $Q_r \cap 2^{-m} \Z^\N$ is a finite set, hence we can write

\[ Q_r \cap 2^{-m} \Z^\N = \{ y_1, ..., y_N \} , \]

where $N\in\N$ is some number depending on both $r$ and $m$. For every $x\in Q_r$ we set

\[ \ca I(x) := \left\{  i \in \{1 , ... , N \} : |x-y_i|_\infty = \min\limits_{1 \leq j \leq N} |x-y_j|_\infty  \right\} . \]

Furthermore, we define

\[ \pi_m^{(r)}(x) := y_{\min \ca I(x)}  .  \]

Observe that the map $\pi_m^{(r)}$ fulfills all the required properties.

\qed

\end{prof}

\begin{de}[Dyadic point]
\xlabel{DE-DYADIC}

We set

\[ \D := \left\{ \left. (x_n)_{n\in\N} \in \R^\N \right| \forall n\in\N, \ \exists m_n \in \N, \ x_n \in 2^{-m_n} \Z^\N \right\} .  \]

We say that $x\in\R^\N$ is a dyadic point if $x \in \D$.

\end{de}

\section{Regularization by Noise}

In this section we are going to prove various estimates regarding the map $\varphi_{n,k}$ defined below.
Surprisingly, although we do not assume any regularity on $b$, $\varphi_{n,k}$ is ``close to'' being Lipschitz continuous in space.
This is due to the noise, which improves the situation significantly. From this point onwards let $(Z^A_t)_{t\in[0,\infty[}$ be an Ornstein--Uhlenbeck process on a probability space $(\Omega, \ca F, (\ca G_t)_{t\in[0,\infty[} , \P)$ with drift term $A$ and initial sigma-algebra $(\ca G_t)_{t\in[0,\infty[}$ as defined in the introduction.

\begin{de}
\xlabel{DE-VARPHI}

Let $b \colon [0,1] \times H \longrightarrow H$ be a Borel measurable function. For $n\in\N$, $k\in\{0, ..., 2^n-1\}$ and $x\in H$ we define
\[ \varphi_{n,k} \colon H \times \Omega \longrightarrow H \]
by
\[ \varphi_{n,k}(b ; x, \omega) := \int\limits_{k2^{-n}}^{(k+1)2^{-n}} b(s, Z^A_s(\omega) + x) - b(s, Z^A_s(\omega)) \d s  . \]

Usually we drop the $b$ and $\omega$ and just write $\varphi_{n,k}(x)$ instead of $\varphi_{n,k}(b ; x, \omega)$. Additionally, we set

\[ \varphi_{n,k}(x,y) := \int\limits_{k2^{-n}}^{(k+1)2^{-n}} b(s, Z^A_s + x) - b(s, Z^A_s + y) \d s  . \]

\end{de}

\begin{re}
\xlabel{RE-METRIC}

Note that for fixed $n\in\N$, $k\in\{0, ..., 2^n-1\}$ and $\omega\in\Omega$ the map

\[ |\varphi_{n,k}(\cdot,\cdot)|_H \colon H \times H \longrightarrow \R_+, \qquad (x,y) \longmapsto |\varphi_{n,k}(x,y)|_H  \]

is a pseudometric on $H$.

\end{re}

\begin{lem}
\xlabel{LOG-CONVEX}

For $r,m \in \N$ we have

\[ \ln(r+m+1)^{1/\gamma} \leq \ln(r+1)^{1/\gamma}  +  \ln(m+1)^{1/\gamma} , \]

where $\gamma$ is the constant from Assumption \ref{ASS}.

\end{lem}

\begin{prof}

Let $r,m \in \N$. We have

\[ r + m +1  \leq rm + r + m + 1 = (r+1)\cdot (m+1) , \]

which implies that

\[ \ln(r + m+1) \leq \ln((r+1)\cdot(m+1)) = \ln(r+1) + \ln(m+1) . \]

Since $\frac1\gamma \leq 1$ we immediately obtain

\[ \ln(r + m+1)^{1/\gamma} \leq \ln(r+1)^{1/\gamma} + \ln(m+1)^{1/\gamma} \]

due to the fact that $x \longmapsto x^{1/\gamma}$ is concave which completes the proof.

\qed

\end{prof}

\begin{thm}
\xlabel{THM-SIGMA}

For every $0 < \eps < \frac16$ there exists $C_\eps \in \R$ such that for every Borel measurable function $b\colon [0,1] \times H \longrightarrow H$ satisfying Assumption \ref{ASS}, $n\in\N \setminus\{0\}$ and $k \in \{ 0 , ..., 2^n - 1 \}$ there exists a measurable set $A_{\eps,b,n,k} \in \ca G_{(k+1)2^{-n}} \subseteq \Omega$ with $\P[A_{\eps,b,n,k}] \leq \frac\eps3 e^{-n}$ such that on $A_{\eps,b,n,k}^c$

\[ |\varphi_{n,k}(x)|_H \leq C_\eps n^{\frac12 + \frac1\gamma} 2^{-n/2} \left( |x|_\infty + 2^{-2^n} \right) \]

holds for all points $x \in 2Q \cap \D$.

\end{thm}

\begin{re}

Note that the constant $C_\eps$ depends on $\eps$ and $\gamma$ from Assumption \ref{ASS}, but not on $b$. Conversely, the set of ``good omegas'' $A_{\eps,b,n,k}^c$ depends on $\eps$, $b$, $n$ and $k$.

\end{re}

\begin{prof}

\textbf{Step 1:}\\

Let $0 < \eps < \frac1{6}$. For $r\geq 0$ recall that $Q_r := \{ x\in Q : |x|_\infty \leq 2 \cdot 2^{-r} \}$. Let $m$ be an integer with $m \geq r$ and $x,y \in 2Q_r \cap 2^{-m} \Z^\N$. We are going to estimate the probability of the event $\{|\varphi_{n,k}(x,y) |_H > \eta\}$ for a suitable $\eta \geq 0$. To this end let $\beta_A > 0$ be the constant from \cite[Corollary 3.1]{Wre16} and for $0 < \eps < \frac1{6}$ we set

\customlabel{SIGMA-ETA-EPS}
\[ \eta_\eps := \sqrt{ \ln\left( \frac{6}\eps \right) } \geq 1  .   \tag{\ref{SIGMA-ETA-EPS}}  \]

Let us consider the following probability.

\[ \P \left[ |\varphi_{n,k}(x,y)|_H  >  \beta_A^{-1/2} \eta_\eps ( 1 + \sqrt{2n}  +  \sqrt{ 9 (1+m-r) } ) \op{ed}(2Q_r)_m |x-y|_\infty 2^{-n/2} \right] . \]

Since $x,y \in 2Q_r \cap 2^{-m} \Z^\N$ and $|\cdot|_\infty$, $|\cdot|_2$ are effectively equivalent norms i.e.~$|\cdot|_2 \leq \sqrt{\op{ed}(2Q_r)_m} |\cdot|_\infty$ (see Proposition \ref{EQUIVALENCE}) the above expression is smaller than

\[ \P \left[ |\varphi_{n,k}(x,y)|_H  >  \beta_A^{-1/2} \eta_\eps ( 1 + \sqrt{2n}  +  \sqrt{ 9 (1+m-r) } ) \sqrt{\op{ed}(2Q_r)_m} |x-y|_2 2^{-n/2} \right] .  \]

Due to Corollary \cite[Corollary 3.1]{Wre16} this probability is smaller than

\[ e^{- \eta_\eps^2 \op{ed}(2Q_r)_m} e^{- \eta_\eps\left( \sqrt{2n}  +  \sqrt{9 (1+m-r)} \right)^2 \op{ed}(2Q_r)_m}  .  \]

Using that $\eta_\eps \geq 1$ and $\op{ed}(2Q_r)_m\geq 1$ the above is bounded from above by

\[ e^{- \eta_\eps^2} e^{- \left( \sqrt{2n}  +  \sqrt{9 (1+m-r)} \right)^2 \op{ed}(2Q_r)_m}   \leq e^{- \eta_\eps^2} e^{- (2n  +  9 (1+m-r)) \op{ed}(2Q_r)_m }  =    e^{- \eta_\eps^2} e^{-2n} e^{- 9 (1+m-r)  \op{ed}(2Q_r)_m }  . \]

In order to get a uniform bound we calculate

\[  \P \left[ \bigcup\limits_{r=0}^{2^n} \bigcup\limits_{m=r}^\infty \!\!\! \bigcup\limits_{\bfrac{x,y \in}{2Q_r \cap 2^{-m} \Z^\N}} \!\!\!\!\!\! \left\{ |\varphi_{n,k}(x,y)|_H > \beta_A^{-1/2} \eta_\eps ( 1 + \sqrt{2n}  +  \sqrt{9 (1+m-r) } ) \op{ed}(2Q_r)_m |x-y|_\infty 2^{-n/2} \right\} \right] \]

\[ \hspace{-27mm}  \leq \sum\limits_{r=0}^{2^n} \sum\limits_{m=r}^\infty \sum\limits_{\bfrac{x,y \in}{2Q_r \cap 2^{-m} \Z^\N}}  e^{- \eta_\eps^2} e^{- 2n} e^{- 9 (1+m-r) \op{ed}(2Q_r)_m} \]

\[ \hspace{10mm}  = e^{- \eta_\eps^2} \sum\limits_{r=0}^{2^n} \sum\limits_{m=r}^\infty \# \{ (x,y) | x,y \in 2Q_r \cap 2^{-m} \Z^\N \}     e^{- 2n} e^{ - 9 (1+m-r) \op{ed}(2Q_r)_m } .  \]

Invoking Theorem \ref{KOLTIK} results in
\[ \# \{ x | x \in 2Q_r \cap 2^{-m} \Z^\N \} \leq \exp \left( 8 ( 1 + m - r ) \op{ed}(2Q_r)_m \right) .  \]

Hence, we can bound the above probability by

\[  e^{- \eta_\eps^2} \sum\limits_{r=0}^{2^n} \sum\limits_{m=r}^\infty   \exp \left( 8 ( 1 + m - r ) \op{ed}(2Q_r)_m \right)  e^{- 2n} e^{- 9(1+m-r) \op{ed}(Q_r)_m } \]

\[ \hspace{-10mm}  = e^{- \eta_\eps^2}  e^{- 2n } \sum\limits_{r=0}^{2^n} \sum\limits_{m=r}^\infty  \exp \left( - ( 1 + m - r ) \op{ed}(2Q_r)_m \right) .  \]

Note that the last sum converges since $\op{ed}(2Q_r)_m \geq 1$ and because of

\[ \sum\limits_{m=r}^\infty  \exp \left( - ( 1 + m - r ) \op{ed}(2Q_r)_m \right)  \leq \sum\limits_{m=0}^\infty  \exp \left( - ( 1 + m ) \right) \leq 1  \]

the above is smaller than

\[  e^{- \eta_\eps^2} \sum\limits_{r=0}^{2^n} e^{-2n}  =  e^{- \eta_\eps^2} (2^n+1)  e^{-2n}  \leq 2 e^{- \eta_\eps^2} e^{-n} . \]

Plugging in Definition \eqref{SIGMA-ETA-EPS} of $\eta_\eps$ the above is smaller than $\frac\eps3 e^{-n}$. In conclusion there exists a measurable set $A_{\eps,b,n,k} \subseteq \Omega$ with $\P[A_{\eps,b,n,k}] \leq \frac\eps3 e^{-n}$ such that on $A_{\eps,b,n,k}^c$ we have

\customlabel{SIGMA-STEP-1}
\begin{align*}\tag{\ref{SIGMA-STEP-1}}
\begin{split}
|\varphi_{n,k}(x,y)|_H  &\leq  \beta_A^{-1/2} \eta_\eps ( 1 + \sqrt{2n}  +  \sqrt{ 9(1+m-r) } ) \op{ed}(2Q_r)_m |x-y|_\infty 2^{-n/2} \\
&\leq  6 \beta_A^{-1/2} \eta_\eps ( \sqrt{n}  +  \sqrt{ 1+m-r } ) \op{ed}(2Q_r)_m |x-y|_\infty 2^{-n/2}
\end{split}
\end{align*}

for $n \geq 1$, $k \in \{ 0, ..., 2^n-1 \}$, $r \in \{ 0, ..., 2^n$\}, $m\geq r$ and $x,y \in 2Q_r \cap 2^{-m} \Z^\N$.\\

\textbf{Step 2:}\\

\textbf{Claim:} For every dyadic number $x \in 2Q_r$ with $r \in \{ 0 , ...,  2^n \}$ and $n \geq 1$, $k \in \{ 0, ..., 2^n-1 \}$ we have

\customlabel{SIGMA-STEP-2}
\[ |\varphi_{n,k}(x)|_H \leq 384 \beta_A^{-1/2} \eta_\eps 2^{-n/2} 2^{-r} \sqrt{n}  (\ln (r+3))^{1/\gamma} .   \tag{\ref{SIGMA-STEP-2}}  \]

on $A_{\eps,b,n,k}^c$. Indeed, let $x$ be a dyadic number such that $x \in 2Q_r$ with $r \in \{ 0 , ...,  2^n \}$. Recall Corollary \ref{RE-PI}. For every $m\in\N$ with $m\geq r$ we set

\[ x_m := 2\pi_{m+1}^{(r)}\left( \frac x2 \right) \in 2Q_r \cap 2^{-m} \Z^\N  , \]

where $\pi_m^{(r)}$ is the map from Corollary \ref{RE-PI}. I.e.~$|x - x_m|_\infty \leq 2^{-m}$. By the triangle inequality (see Remark \ref{RE-METRIC}) and $\varphi_{n,k}(x)=\varphi_{n,k}(x,0)$ we immediately get

\[ |\varphi_{n,k}(x)|_H \leq |\varphi_{n,k}(x_r, 0)|_H  +  \sum\limits_{m=r}^\infty |\varphi_{n,k}(x_{m+1}, x_m)|_H .  \]

Note that the sum on the right-hand side is actually a finite sum, because $x$ is dyadic, so that $x_m = x$ for $m$ sufficiently large. Note that $x_m, x_{m+1} \in 2^{-(m+1)} \Z^\N$ hence, by using inequality \eqref{SIGMA-STEP-1}, the above expression is bounded from above by

\[ 6 \beta_A^{-1/2} \eta_\eps \left( \sqrt{n}  +  \sqrt{ 1+r-r } \right)  \op{ed}(2Q_r)_r |x_r|_\infty 2^{-n/2}   \]
\[  +   6 \beta_A^{-1/2} \eta_\eps \sum\limits_{m=r}^\infty ( \sqrt{n}  +  \sqrt{ 1+(m+1)-r } )   \op{ed}(2Q_r)_{m+1} |x_{m+1} - x_m |_\infty 2^{-n/2} .  \]

Using the definition of $x_m$ and $|x_{m+1} - x_m |_\infty \leq |x_{m+1} - x |_\infty + |x_m - x|_\infty \leq 2^{-m+1}$ this can be estimated from above by

\[ \hspace{-23mm}  6 \beta_A^{-1/2} \eta_\eps \left( \sqrt{n}  +  1 \right)  \op{ed}(2Q_r)_r 2^{-r} 2^{-n/2}   \]
\[ \hspace{19mm}  +   24 \beta_A^{-1/2} \eta_\eps \sum\limits_{m=r}^\infty ( \sqrt{n}  +  \sqrt{ 1+(m+1)-r } ) 2^{-(m+1)} \op{ed}(2Q_r)_{m+1}  2^{-n/2} \]

\[ \leq 24 \beta_A^{-1/2} \eta_\eps 2^{-n/2} \sum\limits_{m=r}^\infty  \left( \sqrt{n}  +  \sqrt{ 1+m-r } \right) \op{ed}(2Q_r)_m 2^{-m} .  \]

Invoking Lemma \ref{EFFDIM} yields that $\op{ed}(Q_r)_m \leq (\ln(m+1))^{1/\gamma}$, where $\gamma>1$ is the constant from Assumption \ref{ASS}. Using this we can further estimate the above expression by

\[ \hspace{11mm}  24 \beta_A^{-1/2} \eta_\eps 2^{-n/2} \sum\limits_{m=r}^\infty   \left( \sqrt{n}  +  \sqrt{1+m-r} \right) (\ln(m+1))^{1/\gamma} 2^{-m}  \]

\[ \hspace{12mm}  \leq 24 \beta_A^{-1/2} \eta_\eps 2^{-n/2} \sum\limits_{m=0}^\infty   \left( \sqrt{n}  +  \sqrt{1+m} \right) (\ln(r+m+1))^{1/\gamma} 2^{-m-r} . \]

Using Lemma \ref{LOG-CONVEX} the above is smaller than

\[ \hspace{-7mm}   24 \beta_A^{-1/2} \eta_\eps 2^{-n/2} 2^{-r} \sum\limits_{m=0}^\infty   \left( \sqrt{n}  +  \sqrt{1+m} \right) \left( (\ln (r+1) )^{1/\gamma} + \ln(m+1)^{1/\gamma} \right) 2^{-m}  \]

\[ \hspace{-15mm}  \leq 24 \beta_A^{-1/2} \eta_\eps 2^{-n/2} 2^{-r} \left[  \sqrt{n} (\ln (r+1))^{1/\gamma} \sum\limits_{m=0}^\infty 2^{-m}    +   \sqrt{n} \sum\limits_{m=0}^\infty (\ln (m+1))^{1/\gamma} 2^{-m}  \right.  \]

\[ \hspace{24mm}  \qquad\qquad \left.  +  (\ln r)^{1/\gamma} \sum\limits_{m=0}^\infty \sqrt{1+m} 2^{-m}  +  \sum\limits_{m=0}^\infty \sqrt{1+m} (\ln (m+1))^{1/\gamma} 2^{-m}  \right] .  \]

Since $\gamma \geq 1$ we can estimate $(\ln (m+1))^{1/\gamma} \leq 2^{m/2}$. The above expression is therefore bounded by

\[ 24 \beta_A^{-1/2} \eta_\eps 2^{-n/2} 2^{-r} \left[  2 \sqrt{n} (\ln (r+1))^{1/\gamma}   +  \sqrt{n} \sum\limits_{m=0}^\infty 2^{-m/2}  +  4 (\ln r)^{1/\gamma}  +  \sum\limits_{m=0}^\infty \sqrt{1+m} 2^{-m/2}  \right]   \]

\[ \leq 24 \beta_A^{-1/2} \eta_\eps 2^{-n/2} 2^{-r} \left[  2 \sqrt{n} (\ln (r+1))^{1/\gamma} + 4 \sqrt{n} + 4 (\ln r)^{1/\gamma} + 6 \right] .  \]

And since we have $1 \leq (\ln(r+3))^{1/\gamma}$, we obtain

\[ |\varphi_{n,k}(x)|_H \leq \underbrace{384 \beta_A^{-1/2} \eta_\eps}_{=: C_\eps} 2^{-n/2} 2^{-r}  \sqrt{n}  (\ln (r+3))^{1/\gamma}  , \]

which proves Claim \eqref{SIGMA-STEP-2}.\\

\textbf{Step 3:}\\

For a fixed $n\in\N$ let $x\in 2Q \cap \D$ such that $|x|_\infty > 2^{-2^{n}}$. We set
\[ r := \lfloor \log_2 |x|_\infty^{-1} \rfloor  \leq  \lfloor 2^{n} \rfloor \leq 2^{n}  .  \]

And hence we have
\[ 2^{-r} =  2^{- \log_2 \lfloor |x|_\infty^{-1} \rfloor} \leq 2^{- \log_2 |x|_\infty^{-1}  + 1}  = 2 |x|_\infty  . \]

Additionally, we have $r\in\{-2, ..., 2^n\}$ and $x \in 2Q_r$, because of the fact that
\[ |x|_\infty  =  2^{- \log_2 |x|_\infty^{-1}}  \leq 2^{-r} . \]

Hence, we can apply Claim \eqref{SIGMA-STEP-2} of Step 2 to obtain

\[ \hspace{-53mm}  |\varphi_{n,k}(x)|_H  \leq  C_\eps  2^{-r} \sqrt{n} 2^{-n/2}  (\ln (r+3))^{1/\gamma}  \]

\[ \hspace{18mm}  \leq C_\eps \sqrt{n} 2^{-n/2}  |x|_\infty \left(\log_2\left(2^{3n}\right)\right)^{1/\gamma} \leq C_\eps \sqrt{n} (3n)^{1/\gamma} 2^{-n/2} |x|_\infty .  \]

\textbf{Step 4:}\\

Conversely to Step 3, for fixed $n\in\N$ let $x \in 2Q \cap \D$ such that $|x|_\infty \leq 2^{-2^{n}}$. Then $x \in Q_r$ with $r = 2^{n}$ so that by Invoking Step 2 (i.e.~inequality \eqref{SIGMA-STEP-2}) we have

\[ \hspace{-10mm}  |\varphi_{n,k}(x)|_H  \leq  C_\eps  2^{-r} \sqrt{n} 2^{-n/2}  (\ln (r+3))^{1/\gamma} \leq  C_\eps  2^{-2^n} \sqrt{n} 2^{-n/2}  \left(\log_2 \left(2^{3n}\right)\right)^{1/\gamma}  \]

\[ \hspace{51.5mm}  \leq C_\eps  \sqrt{n} 2^{-n/2}  2^{-2^{n}}  (3n)^{1/\gamma} .  \]

This concludes the proof.

\qed

\end{prof}

\begin{thm}
\xlabel{THM-RHO}

For every $0 < \eps < 1$ there exists $C_\eps \in \R$ such that for every Borel measurable function $b\colon [0,1] \times H \longrightarrow H$ satisfying Assumption \ref{ASS} there exists a measurable set $A_{\eps,b} \subseteq \Omega$ with $\P[A_{\eps,b}] \leq \eps$ such that on $A_{\eps,b}^c$

\[ |\varphi_{n,k}(x, y)|_H  \leq C_\eps \left[ \sqrt{n} 2^{-n/6} |x - y|_\infty  +  2^{- 2^{\theta n}} \right] \]

holds for all points $x,y \in Q \cap \D$ with $|x-y|_\infty \leq 1$, $n \geq 1$, $k \in \{ 0 , ..., 2^n - 1 \}$ and $\theta := \frac23 \frac\gamma{\gamma+2}$.

\end{thm}

\begin{re}

Note that the constant $C_\eps$ depends on $\eps$ and $\gamma$, but not on $b$. Conversely, the set of ``good omegas'' $A_{\eps,b}^c$ depends on both, $\eps$ and $b$.

\end{re}

\begin{prof}

\textbf{Step 1:}\\

Let $m\in\N$ and $x,y\in Q \cap 2^{-m} \Z^\N$. For $0 < \eps < 1$ we set

\[ \eta_\eps := \sqrt{ \ln\left( \frac1\eps \right) }  \leq  1 .    \]

Analogously to the previous proof we estimate

\[ \P \left[ |\varphi_{n,k}(x, y)|_H  >  \beta_A^{-1/2} \eta_\eps ( 1 + \sqrt{2n}  +  \sqrt{ 5(1+m) } ) \op{ed}(Q)_m |x-y|_\infty 2^{-n/2} \right] ,  \]

where $\beta_A>0$ is the constant from \cite[Corollary 3.1]{Wre16}.
Since $x,y \in Q \cap 2^{-m} \Z^\N$ and $|\cdot|_\infty$, $|\cdot|_2$ are effectively equivalent norms i.e.~\mbox{$|\cdot|_2 \leq \sqrt{\op{ed}(Q)_m} |\cdot|_\infty$} (see Proposition \ref{EQUIVALENCE}) the above expression is smaller than

\[ \P \left[ |\varphi_{n,k}(x, y)|_H  >  \beta_A^{-1/2} \eta_\eps ( 1 + \sqrt{2n}  +  \sqrt{ 5(1+m) } ) \sqrt{\op{ed}(Q)_m} |x-y|_2 2^{-n/2} \right] .  \]

Due to \cite[Corollary 3.1]{Wre16} this expression is bounded by

\[ e^{- \eta_\eps^2 \op{ed}(Q)_m} e^{- \eta_\eps^2 \left( \sqrt{2n}  +  \sqrt{5(1+m)} \right)^2 \op{ed}(Q)_m}  \]

ans since $\eta_\eps \geq 1$ as well as $\op{ed}(Q)_m \geq 1$ the above expression can be estimated from above by

\[ e^{- \eta_\eps^2} e^{- (2n  +  5(1+m)) \op{ed}(Q)_m } \leq e^{- \eta_\eps^2} e^{-2n} e^{- 5(1+m)  \op{ed}(Q)_m }  . \]

Using this, we estimate the following probability

\[  \P \left[ \bigcup\limits_{n=1}^\infty \bigcup\limits_{m=0}^\infty \!\! \bigcup\limits_{\bfrac{x,y\in}{Q \cap 2^{-m} \Z^\N}} \!\! \bigcup\limits_{k=0}^{2^n-1} |\varphi_{n,k}(x, y)|_H > \beta_A^{-1/2} \eta_\eps \left( 1 + \sqrt{2n}  +  \sqrt{5(1+m)} \right) \op{ed}(Q)_m |x-y|_\infty 2^{-n/2} \right] \]

\[ \hspace{-29mm}  \leq \sum\limits_{n=1}^\infty \sum\limits_{m=0}^\infty \sum\limits_{\bfrac{x,y\in}{Q \cap 2^{-m} \Z^\N}} \sum\limits_{k=0}^{2^n-1}  e^{- \eta_\eps^2} e^{- 2n} e^{- 5(1+m) \op{ed}(Q)_m} \]

\[ \hspace{5mm}   \leq e^{- \eta_\eps^2} \sum\limits_{n=1}^\infty \sum\limits_{m=0}^\infty \# \{ (x,y) | x,y \in Q \cap 2^{-m} \Z^\N \}   2^n  e^{- 2n} e^{ - 5(1+m) \op{ed}(Q)_m } .  \]

By invoking Theorem \ref{KOLTIK} for $r=0$ we have
\[ \# \{ (x,y) | x,y \in Q \cap 2^{-m} \Z^\N \} \leq \exp \left( 4 ( 1 + m ) \op{ed}(Q)_m \right) .  \]

So that we can bound the above probability by

\[  \hspace{5mm}  e^{- \eta_\eps^2}  \sum\limits_{n=1}^\infty \sum\limits_{m=0}^\infty  \exp \left( 4 ( 1 + m ) \op{ed}(Q)_m \right)  2^n  e^{- 2n} e^{-  5(1+m) \op{ed}(Q)_m } \]

\[ \hspace{-21.5mm} \leq e^{- \eta_\eps^2} \sum\limits_{n=1}^\infty \sum\limits_{m=0}^\infty 2^n  e^{- 2n }  \exp \left( - ( 1 + m ) \op{ed}(Q)_m \right) .   \]

Note that the last sum converges since $\op{ed}(Q)_m \geq 1$. Hence, the above is bounded from above by

\[ e^{- \eta_\eps^2} \sum\limits_{n=1}^\infty 2^n  e^{- 2n } \underbrace{\sum\limits_{m=0}^\infty  \exp \left( - ( 1 + m ) \right)}_{\leq 1}  \]

so that, in conclusion, we have estimated the above probability by

\[ e^{- \eta_\eps^2} \sum\limits_{n=1}^\infty 2^n  e^{- 2n } \leq e^{- \eta_\eps^2} = \eps .  \]

Therefore, we obtain

\customlabel{RHO-STEP-1}
\begin{align*}\tag{\ref{RHO-STEP-1}}
\begin{split}
|\varphi_{n,k}(x, y)|_H  &\leq  \beta_A^{-1/2} \eta_\eps \left( 1 + \sqrt{2n}  +  \sqrt{5(1+m)} \right) \op{ed}(Q)_m |x-y|_\infty 2^{-n/2} \\
&\leq  5 \beta_A^{-1/2} \eta_\eps \left( \sqrt{n}  +  \sqrt{1+m} \right) \op{ed}(Q)_m |x-y|_\infty 2^{-n/2}
\end{split}
\end{align*}

for $n \geq 1$, $k \in \{ 0, ..., 2^n-1 \}$, $m \in \N$ and for all $x,y \in Q \cap 2^{-m} \Z^\N$ on a set $A_{\eps,b}^c \subseteq \Omega$ with $\P[A_{\eps,b}] \leq \eps$.\\

\textbf{Step 2:}\\

\textbf{Claim:} For all points $x,y \in Q\cap\D$, with $|x-y|_\infty \leq 1$, $n \geq 1$ and $k \in \{ 0, ..., 2^n-1 \}$ we have

\customlabel{RHO-STEP-2}
\[ |\varphi_{n,k}(x,y)|_H \leq 7200 \beta_A^{-1/2} \eta_\eps \sqrt n \left[ 2^{-n/6} |x - y|_\infty  +  2^{- 2^{\theta n}} \right] .   \tag{\ref{RHO-STEP-2}}  \]

on $A_{\eps,b}^c$. Indeed, let $x,y \in Q$ be two dyadic points in $Q$ with $|x-y|_\infty \leq 1$. W.l.o.g.~we assume $x \neq y$.
Fix $m\in\N$ be so that $2^{-m-1} \leq |x - y|_\infty \leq 2^{-m}$. Note that this implies that $m \geq 0$.
Using Corollary \ref{RE-PI} for every $r\in\N$ with $r\geq m$ we set

\[ x_r := \pi_r^{(0)}(x) \in Q \cap 2^{-r} \Z^\N  , \]
\[ y_r := \pi_r^{(0)}(y) \in Q \cap 2^{-r} \Z^\N  . \]

By the triangle inequality (see Remark \ref{RE-METRIC}) we immediately get

\[ |\varphi_{n,k}(x,y)|_H  \leq  |\varphi_{n,k}(x_m, y_m)|_H  +  \sum\limits_{r=m}^\infty |\varphi_{n,k}(x_{r+1}, x_r)|_H  +  \sum\limits_{r=m}^\infty |\varphi_{n,k}(y_{r+1}, y_r)|_H  .  \]

Note that both sums on the right-hand side are actually a finite sums, because $x$ and $y$ are dyadic points. Also note that $x_r, x_{r+1}, y_r, y_{r+1} \in 2^{-(r+1)} \Z^\N$, so that by using inequality \eqref{RHO-STEP-1} the above expression is bounded from above by

\[   5 \beta_A^{-1/2} \eta_\eps \left( \sqrt{n}  +  \sqrt{1+m} \right)  \op{ed}(Q)_{m} |x_m - y_m|_\infty 2^{-n/2}   \]
\[ \hspace{-7.5mm} +   10 \beta_A^{-1/2} \eta_\eps \sum\limits_{r=m}^\infty ( \sqrt{n}  +  \sqrt{ r+2 } )  \op{ed}(Q)_{r+1}  2^{-(r-1)} 2^{-n/2} , \]

where we have used that by the definition of $x_r$ we have $|x_{r+1} - x_r |_\infty \leq |x_{r+1} - x|_\infty + |x_r - x|_\infty \leq 2^{-(r-1)}$ and an analogous calculation for $|y_{r+1} - y_r |_\infty$. Since $|x_m-y_m|_\infty \leq |x_m-x|_\infty + |x-y|_\infty + |y-y_m|_\infty \leq 2^{-(m-2)}$ this can be further estimated from above by

\[ 40 \beta_A^{-1/2} \eta_\eps \sum\limits_{r=m}^\infty ( \sqrt{n}  +  \sqrt{r+1} )  \op{ed}(Q)_{r} 2^{-r} 2^{-n/2} .  \]

Invoking Lemma \ref{EFFDIM} yields that $\op{ed}(Q)_r \leq (\ln(r+1))^{1/\gamma}$, where $\gamma>0$ is the constant from Assumption \ref{ASS}. Using this we can further estimate the above expression by

\[ \hspace{8mm}  40 \beta_A^{-1/2} \eta_\eps \sum\limits_{r=m}^\infty ( \sqrt{n}  +  \sqrt{r+1} )  (\ln (r+1))^{1/\gamma} 2^{-r} 2^{-n/2} \]

and since $\sqrt{n}  +  \sqrt{r+1} \leq 2 \sqrt{n(r+1)}$ this is bounded by

\[ 80 \beta_A^{-1/2} \eta_\eps \sum\limits_{r=m}^\infty \sqrt{n} \sqrt{r+1} (\ln (r+1))^{1/\gamma} 2^{-r} 2^{-n/2}  .  \]

By performing an index shift this can be written as

\[ 80 \beta_A^{-1/2} \eta_\eps \sqrt{n} 2^{-n/2} 2^{-m} \sum\limits_{r=0}^\infty \sqrt{r+m+1} (\ln (r+m+1))^{1/\gamma} 2^{-r} .  \]

We use $\sqrt{r+m+1} \leq \sqrt{r+1} + \sqrt m$ and invoke Lemma \ref{LOG-CONVEX} to estimate this further from above by

\[ 80 \beta_A^{-1/2} \eta_\eps \sqrt{n} 2^{-n/2} 2^{-m} \sum\limits_{r=0}^\infty (\sqrt{r+1} + \sqrt{m}) \left( (\ln(r+1))^{1/\gamma} + (\ln (m+1))^{1/\gamma} \right) 2^{-r} .  \]

Expanding the terms yields

\[ \hspace{-10mm}  80 \beta_A^{-1/2} \eta_\eps \sqrt{n} 2^{-n/2} 2^{-m} \sum\limits_{r=0}^\infty \left[ \sqrt{r+1} (\ln (r+1))^{1/\gamma} + \sqrt{r+1} (\ln (m+1))^{1/\gamma}  \right. \]
\[ \hspace{40mm}  \left.  +  \sqrt{m} (\ln (r+1))^{1/\gamma} +  \sqrt{m} (\ln (m+1))^{1/\gamma} \right] 2^{-r} .  \]

Plugging in $(\ln(m+1))^{1/\gamma} \leq 2^{m/2}$ and evaluating the sum term by term leads us to the following upper bound

\[ 80 \beta_A^{-1/2} \eta_\eps \sqrt{n} 2^{-n/2} 2^{-m}  \left[ \sum\limits_{r=0}^\infty \sqrt{r+1} 2^{-r/2} + (\ln (m+1))^{1/\gamma} \sum\limits_{r=0}^\infty \sqrt{r+1} 2^{-r} \right. \]
\[ \hspace{38mm}  \left. +  \sqrt{m} \sum\limits_{r=0}^\infty 2^{-r/2} +  \sqrt{m} (\ln (m+1))^{1/\gamma} \sum\limits_{r=0}^\infty 2^{-r} \right]   \]

\[ \leq 80 \beta_A^{-1/2} \eta_\eps \sqrt{n} 2^{-n/2} 2^{-m}  \left[ 6 + 3 (\ln (m+1))^{1/\gamma} +  4 \sqrt{m} +  2 \sqrt{m} (\ln (m+1))^{1/\gamma} \right]   \]

\[ \leq 1200 \beta_A^{-1/2} \eta_\eps \sqrt{n} 2^{-n/2}   (\ln (m+1))^{1/\gamma} \sqrt{m+1} 2^{-m} .  \]

\[ \leq 2400 \beta_A^{-1/2} \eta_\eps \sqrt{n} 2^{-n/2}   (m+1)^{1/\gamma} \sqrt{m+1} 2^{-(m+1)} .  \]

In conclusion we finally obtained

\customlabel{RHO-STEP-2-PRE}
\[ |\varphi_{n,k}(x, y)|_H  \leq  2400 \beta_A^{-1/2} \eta_\eps \sqrt{n} 2^{-n/2}   (m+1)^{1/\gamma} \sqrt{m+1} 2^{-(m+1)}  .  \tag{\ref{RHO-STEP-2-PRE}}  \]

We are going to estimate this further using the following claim:\\

Set $\theta := \frac23 \frac\gamma{\gamma+2} > 0$.\\

\textbf{Claim:}\\

For every $n,m \in N$ we have that

\customlabel{RHO-CLAIM}
\[  \sqrt{n} m^{\frac12 + \frac1\gamma} 2^{-m} 2^{-n/2}  \leq  3 \sqrt{n} 2^{-n/6} 2^{-m} + 3\cdot 2^{-2^{\theta n}}    \tag{\ref{RHO-CLAIM}}  \]

holds.

\medskip

\textbf{Proof of Claim \eqref{RHO-CLAIM}:}\\

\textbf{Case 1:} $m \leq 2^{1 + \theta n}$

\[ \sqrt{n} m^{\frac12 + \frac1\gamma} 2^{-n/2} \leq  \sqrt{n} 2^{(1 + \theta n) (\frac12 + \frac1\gamma)} 2^{-n/2} = \sqrt{n} \underbrace{2^{\frac12 + \frac1\gamma}}_{\leq3} \underbrace{2^{\theta (\frac12 + \frac1\gamma) n}}_{= 2^{n/3}} 2^{-n/2} \leq 3 \sqrt{n} 2^{-n/6}  . \]

\textbf{Case 2:} $2^{1 + \theta n} < m$

\[  \underbrace{\sqrt{n} 2^{-n/2}}_{\leq1} m^{\frac12 + \frac1\gamma} 2^{-m} \leq  \underbrace{m^{\frac12 + \frac1\gamma} 2^{-m/2}}_{\leq 3} 2^{-m/2} \leq 3 \cdot 2^{-2^{\theta n}}  . \]

This ends the proof of Claim \eqref{RHO-CLAIM}. Using inequality \eqref{RHO-STEP-2-PRE} and \eqref{RHO-CLAIM} we conclude that

\[ |\varphi_{n,k}(x,y)|_H \leq 7200 \beta_A^{-1/2} \eta_\eps \sqrt n \left[ 2^{-n/6}  2^{-m}  +  2^{- 2^{\theta n}} \right] .  \]

Recall that $2^{-m-1} \leq |x - y|_\infty$ so that the above is smaller than

\[ 14400 \beta_A^{-1/2} \eta_\eps \left[ \sqrt{n} 2^{-n/6} |x - y|_\infty  +  2^{- 2^{\theta n}} \right]  ,  \]

which finishes the proof of Claim \eqref{RHO-STEP-2} and hence the assertion.

\qed

\end{prof}

\section{Continuity of the map $\boldsymbol{\varphi_{n,k}}$}

In this section we will prove that for almost all $\omega\in\Omega$ the map

\[ \varphi_{n,k}(\omega) \colon H \longrightarrow H \]

is continuous. Furthermore, we will show that on a suitable class of Lipschitz functions $h$ and their dyadic piecewise approximations that the map

\[ h \longmapsto \int\limits_0^1 b(t, Z_t(\omega) + h(t)) \d t  \]

is continuous w.r.t.~the maximum norm.

\begin{de}
\xlabel{DE-QA}

We set

\[ Q^A := \{ x = (x_n)_{n\in\N} \in\R^\N | \lambda_n e^{2\lambda_n} x_n^2  \leq  1 \} ,  \]

where $(\lambda_n)_{n\in\N}$ are the eigenvalues of the operator $A$ of our Ornstein--Uhlenbeck process $Z^A$.

\end{de}

\begin{de}
\xlabel{DE-PHI}

We define

\begin{align*}
\Phi   &:= \{ h \colon [0,1] \longrightarrow Q \cap Q^A \colon |h(s) - h(t)|_\infty \leq 2|s-t|,\ \forall s,t\in[0,1] \} ,  \\
\Phi_n &:= \left\{ h \colon [0,1] \longrightarrow Q \cap Q^A \cap \D \left| \bfrac{\forall 0\leq k < 2^n \colon \forall s,t\in [k2^{-n}, (k+1)2^{-n}[ \colon h(s)=h(t) \text{ and}}{\forall m,\ell\in\Z\cap [0,2^n] \colon |h(m2^{-n}) - h(\ell 2^{-n})|_\infty \leq 2|m-\ell|2^{-n}} \right.\right\} ,  \\
\Phi^* &:= \Phi \cup \bigcup\limits_{n=1}^\infty \Phi_n .
\end{align*}

\end{de}

\begin{re}

Note that elements in $\Phi$ are continuous, since functions in $\Phi$ are Lipschitz continuous (with Lipschitz constant at most $2$). $\Phi_n$ will be used to approximate elements in $\Phi$. Also note that $\Phi$ and $\Phi_n$ are separable w.r.t.~the maximum norm and hence $\Phi^*$ is separable.

\medskip

Observe that the above spaces are constructed in such a way that the assumptions we impose on $b$ (see Assumption \ref{ASS}) implie that the function $u$ from Proposition \ref{PRO-GIRSANOV} is in the space $\Phi$.
I.e.~the difference of two solutions of \eqref{SDE} always lives in the space $\Phi$ due to Assumption \ref{ASS}.

\end{re}

\begin{lem}
\xlabel{APPROX-LEM}

Let $h\in\Phi^*$ and $n\in\N$. We then have

\[ \sum\limits_{k=0}^{2^n-1} \left| h((2k+1)2^{-(n+1)}) - h(2k2^{-(n+1)}) \right|_\infty \leq 1 .  \]

\end{lem}

\begin{prof}

Let $h\in\Phi^*$ and $n\in\N$ be as in the assertion. If $h\in\Phi$ the inequality follows immediately from the Lipschitz continuity of $h$. Let $h\in\Phi_m$ for some $m\in\N$.\\

\textbf{Case 1:} $m\geq n+1$

We have

\[ \sum\limits_{k=0}^{2^n-1} \left| h((2k+1) 2^{-(n+1)}) - h(2k 2^{-(n+1)}) \right|_\infty  \]

\[ =  \sum\limits_{k=0}^{2^n-1} \left| h((2k+1) 2^{m-(n+1)}2^{-m}) - h(2k 2^{m-(n+1)}2^{-m}) \right|_\infty  . \]

Using the assumption that $h\in\Phi_m$ by definition of $\Phi_m$ the above expression us bounded from above by

\[  \sum\limits_{k=0}^{2^n-1} 2 \cdot 2^{m-(n+1)}2^{-m} = 1 .  \]

\textbf{Case 2:} $m < n+1$

Since $h\in\Phi_m$ is constant on all intervals of the form $[k2^{-m} , (k+1)2^{-m}[$ the sum simplifies to

\[ \sum\limits_{k=0}^{2^n-1} \left| h((2k+1) 2^{-(n+1)}) - h(2k 2^{-(n+1)}) \right|_\infty = \sum\limits_{k=0}^{2^{m-1}-1} \left| h((2k+1) 2^{-m}) - h(2k 2^{-m}) \right|_\infty .  \]

And using the definition of $\Phi_m$ the above sum is bounded by

\[ \sum\limits_{k=0}^{2^{m-1}-1} 2\cdot2^{-m} = 1 .  \]

\qed

\end{prof}

\begin{lem}
\xlabel{APPROX}

For every $0 < \eps < 1$ there exist $\delta > 0$ such that for every open set $U \subseteq [0,1] \times H$ with mass $\mu[U] \leq \delta$, where $\mu = \mathrm dt \otimes Z^A_t(\P)$ there is a measurable set $\Omega_{\eps,U} \subseteq \Omega$ with

\[ \P[\Omega \setminus \Omega_{\eps,U}] \leq \eps \]

such that the inequality

\[ \int\limits_0^1 \bbm 1_U( s, Z^A_s + h(s)) \d s \leq \eps \]

holds on $\Omega_{\eps,U}$ uniformly for any $h\in \Phi^*$.

\end{lem}

\begin{prof}

Let $0 < \eps < 1$ and let $C_{\eps/2}$ be the constant from Theorem \ref{THM-RHO}. Recall that $\theta := \frac23 \frac\gamma{\gamma+2}$. Choose $m\in\N$ sufficiently large, so that

\customlabel{APPROX-LEMMA-DEF-M}
\[ 6C_{\eps/2} \sum\limits_{n=m}^\infty \sqrt n 2^{-n/6} \leq \frac\eps2  \qquad \text{and} \qquad  m \geq \frac4{\theta^2 \ln(2)^2} .  \tag{\ref{APPROX-LEMMA-DEF-M}}   \]

Set $\ca N_m := Q \cap Q^A \cap 2^{-m} \Z^\N$ and note that $\ca N_m$ is a finite $2^{-m}$-net of $Q \cap Q^A$ w.r.t.~the maximum norm.

\medskip

Also, observe that $Z^A_t(\P)$ is equivalent to the invariant measure $N(0, \frac12 A^{-1})$ due to \cite[Theorem 11.13]{DZ92} and analogously $(Z_t^{A} + h(t))(\P)$ to $N(h(t), \frac1{2} A^{-1})$. Let $z\in\ca N_m$ then $z$ is in the domain of $A$ because

\[ \sum\limits_{n\in\N} \la z,e_n\ra^2 \lambda_n^2  = \sum\limits_{n\in\N} |z_n|^2 \lambda_n^2 < \infty  . \]

Set $y := 2 A z$ then $y\in H$ hence, \cite[Corollary 2.4.3]{Bog98} is applicable which implies that the measures $N(0, \frac12 A^{-1})$ and $(Z^A+z)(\P)$ are equivalent. We set

\[ \mu := \mathrm dt \otimes Z^A_t(\P), \qquad \mu_z := \mathrm dt \otimes (Z^A_t + z)(\P) \]

for all $z\in D(A)$. By the Radon--Nikodyn Theorem there exist densities $\rho_z$ so that

\[ \frac{\mrm d\mu_z}{\mrm d\mu} = \rho_z . \]

Furthermore, the family $\{ \rho_z | z\in\ca N_m \}$ is uniformly integrable, since $\ca N_m$ is finite. Hence, there exists $\delta > 0$ such that

\customlabel{APPROX-LEMMA-DENSITIES}
\[ \int\limits_A \rho_z(t,x) \d\mu(t,x) \leq \frac{\eps^2}{4 \cdot 2^m \#(\ca N_m)}, \qquad \forall z\in\ca N_m  \tag{\ref{APPROX-LEMMA-DENSITIES}}  \]

for every measurable set $A \subseteq \Omega$ with $\mu[A] \leq \delta$. Let $U \subseteq [0,1]\times H$ be open with mass $\mu[U] \leq \delta$.
Then, by invoking Theorem \ref{THM-RHO} for the function $\bbm 1_U$ with the constant $C_{\eps/2}$, there exists a measurable set $A_{\eps,U} \subseteq \Omega$ with $\P[A_{\eps,U}] \leq \frac\eps2$ such that

\[ \left| \int\limits_{k 2^{-n}}^{(k+1) 2^{-n}} \!\! \bbm 1_U(t, Z^A_t + x)  -  \bbm 1_U(t, Z^A_t + y)  \d t \right|  \leq  C_{\eps/2} \left( \sqrt n 2^{-n/6} |x-y|_\infty  +  2^{-2^{\theta n}} \right)  . \]

holds for every $n\geq 1$, $k\in \{0, ..., 2^n -1 \}$ and $x,y \in Q \cap \D$ on $A_{\eps,U}^c$. Furthermore, we define the events $B_{\eps,U}$ by

\[ B_{\eps,U} := \bigcup\limits_{z\in\ca N_m} \left\{  \int\limits_0^1 \bbm 1_U (s, Z^A_s + z) \d s > \frac\eps{2\cdot 2^m}  \right\} .  \]

We then have

\[ \hspace{-60mm}  \P[B_{\eps,U}] = \P \left[  \bigcup\limits_{z\in\ca N_m} \left\{ \int\limits_0^1 \bbm 1_U (s, Z_s + z) \d s > \frac\eps{2\cdot 2^m} \right\} \right] \]

\[  \hspace{12mm}  \leq \sum\limits_{z\in\ca N_m}  \P \left[  \int\limits_0^1 \bbm 1_U (s, Z^A_s + z) \d s > \frac\eps{2\cdot2^m} \right]  \leq  \frac{2\cdot2^m}\eps \sum\limits_{z\in\ca N_m}  \E \int\limits_0^1 \bbm 1_U (s, Z^A_s + z) \d s \]

\[ \hspace{2mm}  = \frac{2\cdot2^m}\eps \sum\limits_{z\in\ca N_m}  \int\limits_{[0,1] \times H} \bbm 1_U (s, x) \d \mu_z(s, x) = \frac{2\cdot2^m}\eps \sum\limits_{z\in\ca N_m} \int\limits_U \rho_z(s,x) \d\mu(s, x) . \]

Since $\mu[U] \leq \delta$ using inequality \eqref{APPROX-LEMMA-DENSITIES} the above is bounded from above by

\[ \frac{2\cdot2^m}\eps \#(\ca N_m) \frac{\eps^2}{4\cdot 2^m \#(\ca N_m)} = \frac\eps2 .  \]

In conclusion we proved that we have $\P[ B_{\eps,U} ] \leq \frac\eps2$ and therefore obtained that

\[ \P[ A_{\eps,U}^c \cap B_{\eps,U}^c] \geq 1 - \eps .  \]

For every $h\in\Phi$ and $n\in\N$ we define

\customlabel{APPROX-LEMMA-H_N}
\[ h_n(t) := \sum\limits_{k=0}^{2^n-1}  \bbm 1_{[k2^{-n}, (k+1)2^{-n}[}(t) \frac{\lfloor 2^n h(k2^{-n}) \rfloor}{2^n} \in \underbrace{Q \cap Q^A \cap 2^{-n} \Z^\N}_{=\ca N_n} ,  \qquad \forall t \in [0,1]  ,  \tag{\ref{APPROX-LEMMA-H_N}} \]

where $\lfloor \cdot \rfloor$ denotes the componentwise floor function. Note that $h_n$ is $Q \cap Q^A$-valued since $h$ is $Q \cap Q^A$-valued. Furthermore, $h_n(t)$ is a dyadic number for all $t\in [0,1]$. Also note that $h_n$ converges to $h$ for $n\rightarrow\infty$.

\medskip

Now, let

\[ E_{\eps,U} := \bigcap\limits_{h \in \Phi^*} \left\{ \int\limits_0^1 \bbm 1_U(t, Z^A_t + h(t)) \d t  \leq \eps \right\}  . \]

We are going to prove that $A_{\eps,U}^c \cap B_{\eps,U}^c \subseteq E_{\eps,U}$ holds. To this end let $\omega \in A_{\eps,U}^c \cap B_{\eps,U}^c$. Using that $\omega \in B_{\eps,U}^c$ we have

\[ \left| \int\limits_0^1 \bbm 1_U(t, Z^A_t(\omega) + h_m(t)) \d t \right|  \leq  \sum\limits_{k=0}^{2^m-1} \left| \int\limits_{k2^{-m}}^{(k+1)2^{-m}}  \right. \!\!\! \bbm 1_U(t, Z^A_t(\omega) + \underbrace{h_m(t)}_{\in\ca N_m}) \d t \left. \vphantom{\int\limits_{k2^{-m}}^{(k+1)2^{-m}}} \right|  \leq  \sum\limits_{k=0}^{2^m-1}  \frac{\eps}{2 \cdot 2^m}  =  \frac\eps2 .  \]

And since $\omega\in A_{\eps,U}^c$ we obtain for $n\geq m$

\[ \hspace{-35mm}  \left| \int\limits_0^1 \bbm 1_U(t, Z^A_t(\omega) + h_{n+1}(t)) - \bbm 1_U(t, Z^A_t(\omega) + h_n(t)) \d t \right| \]

\[ \hspace{7mm}  \leq \sum\limits_{k=0}^{2^{n+1}-1} \left| \int\limits_{k2^{-(n+1)}}^{(k+1)2^{-(n+1)}} \!\! \right. \bbm 1_U(t, Z^A_t(\omega) + \underbrace{h_{n+1}(t)}_{\in Q \cap \D}) - \bbm 1_U(t, Z^A_t(\omega) + \underbrace{h_n(t)}_{\in Q \cap \D}) \d t \left. \vphantom{\int\limits_{k2^{-(n+1)}}^{(k+1)2^{-(n+1)}}} \right| \]

\[ \hspace{-5mm}  \leq \sum\limits_{k=0}^{2^{n+1}-1} C_{\eps/2} \left( \sqrt{n} 2^{-n/6} |h_{n+1}(k2^{-n-1}) - h_n(k2^{-n-1})|_\infty + 2^{-2^{\theta n}} \right) \]

\[ \hspace{8mm}  \leq C_{\eps/2} \left[ 2^{n+1} 2^{-2^{\theta n}}  + \sqrt n 2^{-n/6}  \sum\limits_{k=0}^{2^{n+1}-1} |h_{n+1}(k2^{-n-1}) - h_n((k/2)2^{-n})|_\infty   \right] .  \]

Note that since $h_n$ is constant on intervals of the form $[k2^{-n}, (k+1)2^{-n}[$ we have $h_n((k/2)2^{-n}) = h_n( \lfloor k/2 \rfloor 2^{-n})$, so that the above equals

\[ C_{\eps/2}  \left[ 2^{n+1} 2^{-2^{\theta n}} + \sqrt n 2^{-n/6} \sum\limits_{k=0}^{2^{n+1}-1} |h_{n+1}(k2^{-n-1}) - h_n(\lfloor k/2 \rfloor 2^{-n})|_\infty \right]  .  \]

Plugging in Definition \eqref{APPROX-LEMMA-H_N} yields that the above expression can be written as

\[  C_{\eps/2}  \left[ 2^{n+1-2^{\theta n}} + \sqrt n 2^{-n/6} \sum\limits_{k=0}^{2^{n+1}-1} 2^{-n-1} \left| \lfloor 2^{n+1}h(k2^{-n-1}) \rfloor - 2\left\lfloor 2^n h\left( \left\lfloor k/2 \right\rfloor 2^{-n} \right) \right\rfloor \right|_\infty \right] \]

\[  \leq C_{\eps/2} \left[ 2^{n - \frac12 \theta^2 \ln(2)^2 n^2} + \sqrt n 2^{-n/6} \sum\limits_{k=0}^{2^{n+1}-1} \right. 2^{-n-1} \underbrace{\left| \lfloor 2^{n+1}h(k2^{-n-1}) \rfloor - 2^{n+1}h(k2^{-n-1}) \right|_\infty}_{\leq 1} \]

\[  \hspace{13.5mm}  \left. + \sqrt n 2^{-n/6} \sum\limits_{k=0}^{2^{n+1}-1} \left| h(k2^{-n-1}) - h\left( \left\lfloor k/2 \right\rfloor  2^{-n} \right) \right|_\infty \right. \]

\[  \hspace{43.4mm}  + \sqrt n 2^{-n/6} \sum\limits_{k=0}^{2^{n+1}-1} 2^{-n} \underbrace{\left| 2^n h\left( \left\lfloor k/2 \right\rfloor  2^{-n} \right) - \left\lfloor 2^n h\left( \left\lfloor k/2 \right\rfloor  2^{-n} \right) \right\rfloor \right|_\infty}_{\leq1} \left. \vphantom{\sum\limits_{k=0}^{2^{n+1}-1}} \right] \]

\[ \leq C_{\eps/2}  \left[ 2^{n - \frac12 \theta^2 \ln(2)^2 m n} + 3 \sqrt n 2^{-n/6} + \sqrt n 2^{-n/6} \sum\limits_{k=0}^{2^{n+1}-1} \left| h\left(k2^{-(n+1)}\right)  -  h\left( 2\left\lfloor k/2 \right\rfloor  2^{-(n+1)} \right) \right|_\infty \right] .  \]

Since $k = 2 \lfloor k/2 \rfloor$ in case $k$ is even the sum can be restricted to $k$ of the form $k = 2k' + 1$ for $k' \in \{ 0, ..., 2^n -1 \}$. with the help of \eqref{APPROX-LEMMA-DEF-M} the above is bounded by

\[ C_{\eps/2}  \left[ 2^{n - 2n} + 3 \sqrt n 2^{-n/6} + \sqrt n 2^{-n/6} \sum\limits_{k'=0}^{2^n-1} \left| h\left((2k'+1)2^{-(n+1)}\right)  -  h\left( 2k' 2^{-(n+1)} \right) \right|_\infty \right] .  \]

Using Lemma \ref{APPROX-LEM} we can further estimate the above sum by $1$ so that in conclusion we obtain

\[ \left| \int\limits_0^1 \bbm 1_U(t, Z^A_t(\omega) + h_{n+1}(t)) - \bbm 1_U(t, Z^A_t(\omega) + h_n(t)) \d t \right| \leq 6 C_{\eps/2} \sqrt{n} 2^{-n/6}  .  \]

Therefore as long as $\omega \in A_{\eps,U}^c \cap B_{\eps,U}^c$ we have by Lebesgue's dominated convergence Theorem, the lower semi-continuity of $\bbm 1_U$ and by the above calculation

\[ \int\limits_0^1 \bbm 1_U(t, Z^A_t(\omega) + h(t)) \d t \leq \lim\limits_{n\rightarrow\infty} \int\limits_0^1 \bbm 1_U(t, Z^A_t(\omega) + h_n(t)) \d t \]

\[ = \int\limits_0^1 \bbm 1_U(t, Z^A_t(\omega) + h_m(t)) \d t + \sum\limits_{n=m}^\infty \int\limits_0^1 \bbm 1_U(t, Z^A_t(\omega) + h_{n+1}(t)) - \bbm 1_U(t, Z^A_t(\omega) + h_n(t)) \d t \]

\[ \leq \frac\eps2 + 6 C_{\eps/2} \sum\limits_{n=m}^\infty \sqrt n 2^{-n/6}  \overset{\eqref{APPROX-LEMMA-DEF-M}}\leq \frac\eps2 + \frac\eps2 = \eps  . \]

In conclusion we have proven that $A_{\eps,U}^c \cap B_{\eps,U}^c \subseteq E_{\eps,U}$ and hence $\P[E_{\eps,U}] \geq 1 - \eps$ which completes the proof.

\qed

\end{prof}

\begin{thm}[Approximation Theorem]
\xlabel{THM-APPROX}

Let $b \colon [0,1] \times H \longrightarrow H$ be a bounded, Borel measurable function satisfying Assumption \ref{ASS}. There exists a measurable set $\Omega'\subseteq \Omega$ with $\P[\Omega'] = 1$ such that for every sequence $(h_m)_{m\in\N} \subseteq \Phi^*$ which converges pointwise to a function $h\in\Phi^*$ i.e.~$\lim\limits_{m\rightarrow\infty} |h(t)-h_m(t)|_H=0$ we have

\[ \lim\limits_{m\rightarrow\infty}  \int\limits_0^1 b(s, Z^A_s + h_m(s)) \d s  =  \int\limits_0^1 b(s, Z^A_s + h(s)) \d s  \]

on $\Omega'$.

\end{thm}

\begin{prof}

Let $b$ be as in the assertion. For $\ell\in\N$ let $\eps_\ell := 2^{-\ell}$. By Lemma \ref{APPROX} for every $\eps_\ell$ there exists for every $\ell\in\N$ a $\delta_\ell$ such that for every pair $(\eps_\ell, \delta_\ell)$ the conclusions of Lemma \ref{APPROX} holds. Applying Lusin's Theorem to the pair $(b, \delta_\ell)$ yields for every $\ell\in\N$ a closed set $K_\ell \subseteq [0,1]\times H$ with $\mu[K_\ell^c] \leq \delta_\ell$, where $\mu := \mathrm dt \otimes Z^A_t(\P)$, so that

\[ b\mid_{K_\ell} \colon K_\ell \longrightarrow H, \qquad (t,x) \longmapsto b(t,x) \]

is continuous. By Dugundji's Extension Theorem (see \cite[Theorem 4.1]{Dug51}) (applied to the above maps) there exist functions $\bar b_\ell \colon [0,1] \times H \longrightarrow H$ such that

\[ \hspace{25mm}  b (t, x) = \bar b_\ell (t, x ), \qquad \forall (t,x) \in K_\ell , \]
\[ \hspace{-10mm}  \|\bar b_\ell\|_\infty \leq 1 \]

and

\[ \bar b_\ell  \text{ is continuous.} \]

Then, by invoking Lemma \ref{APPROX} for $(\eps_\ell, \delta_\ell, K_\ell^c)$ we obtain for every $\ell\in\N$ a measurable set $\Omega'_\ell$ with $\P[\Omega'_\ell] \geq 1 - \eps_\ell$ such that for any $\omega \in \Omega'_\ell$ and $h\in \Phi^*$

\[ \int\limits_0^1 \bbm 1_{K_\ell^c}(s, Z^A_s(\omega) + h(s)) \d s \leq \eps_\ell \]

holds. Let

\[ \Omega' := \liminf\limits_{\ell\rightarrow\infty} \Omega'_\ell .  \]

Since we have
\[ \sum\limits_{\ell\in\N} \P[ \Omega'^c_\ell ] \leq \sum\limits_{\ell\in\N} \eps_\ell = \sum\limits_{\ell\in\N} 2^{-\ell} < \infty \]

the Borel--Canteli Lemma implies that

\[ \P[ \limsup\limits_{\ell\rightarrow\infty} \Omega'^c_\ell ] = 0  \quad \Rightarrow  \quad \P[\Omega'] = 1 . \]

Let $\omega \in \Omega'$ be fixed. Then, there is an $N(\omega) \in\N$ such that for all $\ell>N(\omega)$ we have $\omega \in \Omega_\ell$ and therefore for all $m\in\N$ we obtain

\customlabel{APPROX-THM-INVOC}
\[ \left| \int\limits_0^1 \bbm 1_{K_\ell^c} (s, Z^A_s(\omega) + h_m(s)) \d s \right| \leq \eps_\ell \tag{\ref{APPROX-THM-INVOC}} .  \]

Note that inequality \eqref{APPROX-THM-INVOC} also holds if we replace $h_m$ by $h$, since $h\in\Phi^*$ by assumption.

The assertion now follows easily by the following calculation

\[ \left| \int\limits_0^1 b(s, Z^A_s(\omega) + h_m(s))   -  \bar b_\ell(s, Z^A_s(\omega) + h_m(s))  \d s \right|_H   \]

\[ \leq \int\limits_0^1 \bbm 1_{K_\ell^c} (s, Z^A_s(\omega) + h_m(s)) \underbrace{\left| b(s, Z^A_s(\omega) + h_m(s)) - \bar b_\ell(s, Z^A_s(\omega) + h_m(s)) \right|_H}_{\leq 2} \d s  \]

\[ \leq 2 \underbrace{\int\limits_0^1 \bbm 1_{K_\ell^c} (s, Z^A_s(\omega) + h_m(s)) \d s}_{\leq \eps_\ell \text{ by \eqref{APPROX-THM-INVOC}}}  .  \]

In conclusion we have

\[ \lim\limits_{m\rightarrow\infty}  \left|  \int\limits_0^1  b(s, Z^A_s(\omega) + h_m(s)) -  b(s, Z^A_s + h(s))  \d s  \right|_H  \]

\[ \leq \lim\limits_{m\rightarrow\infty}  \left|  \int\limits_0^1  b(s, Z^A_s(\omega) + h_m(s))   -  \bar b_\ell(s, Z^A_s + h_m(s))  \d s  \right.  \]

\[ \hspace{13mm}  +  \left. \int\limits_0^1  \bar b_\ell(s, Z^A_s(\omega) + h_m(s))  -   b(s, Z^A_s(\omega) + h(s))  \d s  \right|_H  .  \]

Using the above calculation this is bounded from above by

\[ 2 \eps_\ell + \lim\limits_{m\rightarrow\infty}  \left| \int\limits_0^1 \bar b_\ell(s, Z^A_s(\omega) + h_m(s))  -   b(s, Z^A_s(\omega) + h(s))  \d s \right|_H  .  \]

Since $\bar b_\ell$ is continuous and $h_m$ converges pointwise to $h$ this is the same as

\[ 2 \eps_\ell  +  \left| \int\limits_0^1 \bar b_\ell(s, Z^A_s(\omega) + h(s)) \d s  -  b(s, Z^A_s(\omega) + h(s)) \d s \right|_H  \]

\[ \leq 2 \eps_\ell + \int\limits_0^1 \bbm 1_{K_\ell^c} (s, Z^A_s(\omega) + h(s)) \underbrace{\left| \bar b_\ell(s, Z^A_s(\omega) + h(s)) - b(s, Z^A_s(\omega) + h(s)) \right|_H}_{\leq 2} \! \d s  \leq  4 \eps_\ell ,  \]

where the last inequality follows by invoking inequality \eqref{APPROX-THM-INVOC} for $h_m$ replaced by $h$. Taking the limit $\ell\rightarrow\infty$ completes the proof of the assertion, since the left-hand side is independent of $\ell$.

\qed

\end{prof}

Using the above Approximation Theorem we can now extend the estimates obtained in Section 2 to the whole space $Q$ as the following Corollary shows.

\begin{cor}
\label{COR-SIGMA-RHO}

For every $0 < \eps < \frac16$ there exists $C_\eps \in \R$ such that for every function $b\colon [0,1] \times H \longrightarrow H$ satisfying Assumption \ref{ASS}, $n\in\N$ and $k \in \{ 0 , ..., 2^n - 1 \}$ there exists a measurable set $A_{\eps,b,n,k} \in \ca G_{(k+1)2^{-n}} \subseteq \Omega$ with $\P[A_{\eps,b,n,k}] \leq \frac\eps3 e^{-n}$ such that

\[ \bbm 1_{A_{\eps,b,n,k}^c} |\varphi_{n,k}(x)|_H \leq C_\eps n^{\frac12 + \frac1\gamma} 2^{-n/2} \left( |x|_\infty + 2^{-2^n} \right) \]

holds for every $x\in 2Q$ and by setting

\[ A_{\eps,b} := \bigcup\limits_{n=1}^\infty \bigcup\limits_{k=0}^{2^n-1} A_{\eps,b,n,k} \]

we have $\P[A_{\eps,b}] \leq \eps$ with the property that

\[ \bbm 1_{A_\eps^c} |\varphi_{n,k}(x, y)|_H  \leq C_\eps \left[ \sqrt{n} 2^{-n/6} |x - y|_\infty  +  2^{- 2^{\theta n}} \right] \]

holds for all $x,y \in 2Q$, $n \geq 1$ and $k \in \{ 0 , ..., 2^n - 1 \}$, where $\theta := \frac23 \frac{\gamma}{\gamma+2}$.

\end{cor}

\begin{prof}

The first inequality follows from Theorem \ref{THM-SIGMA} for all points $x \in 2Q \cap \D$. For general points $x \in 2Q$ this follows by approximating $2Q \cap \D \ni x_n \longrightarrow x$ and using Theorem \ref{THM-APPROX}.

\medskip

The second inequality follows in the same way by combining Theorem \ref{THM-RHO} and Theorem \ref{THM-APPROX}. Note that the estimate can be trivially extended from points $x,y \in Q$ with $|x-y|_\infty \leq 1$ to $x,y \in 2Q$ by changing the constant $C_\eps$ and using that $\varphi_{n,k}$ is a seminorm.

\medskip

Observe that one can choose ($C_\eps$ / $A_{\eps,b}$), so that the conclusion of Theorem \ref{THM-SIGMA} and \ref{THM-RHO} hold (with the same constant / one the same set).

\qed

\end{prof}

\section{Long-time Regularization by Noise via Euler Approximation}

In this section we will prove estimates for terms of the type

\[ \sum\limits_{q=1}^N |\varphi_{n,k+q}(x_{q+1}, x_q)|_H  .  \]

We will first prove a concentration of measure result for the above term in Lemma \ref{LEM-EULER}.
Using this we prove a $\P$-a.s.~sure version of this estimate in Theorem \ref{THM-EULER}.
However, this estimate only holds for medium-sized $N$.
By splitting the sum and using Theorem \ref{THM-EULER} repetitively we conclude the full estimate in Corollary \ref{COR-GLUING}.

\medskip

Note that applying Corollary \ref{COR-SIGMA-RHO} to every term under the sum would result in an estimate of order $\ca O(\sqrt n 2^{-n/6} N)$. Since $N$ will later be chosen to be of order $2^n$ this is of no use. The technique to overcome this is two-fold:

\medskip

On the one hand the $\varphi_{n,k+q}$ terms have to ``work together'' to achieve an expression of order $\ca O(N)$. However, since $\{ \varphi_{n,k+q}(x_q) | q=1, ..., N\}$ are ``sufficiently uncorrelated'' the law of large numbers tells us to expect on average an estimate of order $\ca O(\sqrt N)$.

\medskip

On the other hand in later applications $x_q$ will be values from the solution of the integral equation $\text{(IE)}_{\omega}$, so that it is reasonable to assume that $|x_{q+1} - x_q|_H \approx |\varphi_{n,k+q}(x_q)|_H$. Exploiting this enables to use \textit{both} of our previous established estimates for every $|\varphi_{n,k+q}(x_{q+1},x_q)|_H$ term.

\medskip

Using both techniques we end up with an estimate of order $\ca O(2^{-n} N)$ (see Corollary \ref{COR-GLUING}).

\begin{thm}[Burkholder--Davis--Gundy Inequality]

Let $(M_n, \ca F_n)_{n\in\N}$ be a real-valued martingale. For $2 \leq p < \infty$ we have

\customlabel{BDG-1}
\[ (\E |M_n|^p)^{1/p} \leq p (\E |\la M \ra_n^{p/2} )^{1/p} .  \tag{\ref{BDG-1}}  \]

\end{thm}

\begin{prof}

In the celebrated paper \cite[Section 3]{Dav76} it is shown that the optimal constant in our case is the largest positive root of the Hermite polynomial of order $2k$. We refer to the appendix of \cite{Ose12} for a discussion of the asymptotic of the largest positive root. See also \cite[Appendix B]{Kho14}, where a self-contained proof of the Burkholder--Davis--Gundy Inequality with asymptotically optimal constant can be found for the one-dimensional case.

\end{prof}

\begin{lem}
\xlabel{LEMMA-MARTINGALE-2}

Let $(M_n)_{n\in\N}$ be a martingale of the form

\[ M_r := \sum\limits_{k=1}^r X_k  \]

with $\E [X_k^{p}] \leq C^p p^p$ for all $k\in\N$ and $p\in[1,\infty[$ then

\[ \E\left[ \exp\left( \frac18 \left(\frac{M_r}{C \sqrt r}\right)^{1/2} \right) \right] \leq 2  \]

holds for all $r\in\N$.

\end{lem}

\begin{prof}

Let $(M_n)_{n\in\N}$ be a martingale. Using the Burkholder--Davis--Gundy Inequality \eqref{BDG-1} for every $r,p\in\N$ with $p\geq 2$ we have

\begin{align*}
 \E [ M_r^p ] &\leq p^{p} \E [ \la M \ra_r^{p/2} ] = p^{p} \E \left[ \left( \sum\limits_{k=1}^r X_k^2 \right)^{p/2} \right] \\
 &\leq p^{p} r^{p/2-1} \E \left[ \sum\limits_{k=1}^r X_k^{p} \right] \leq p^{p} r^{p/2-1} r C^p p^p = C^p r^{p/2} p^{2p}  .
\end{align*}

In conclusion we obtain

\customlabel{MART2-3}
\[ \E [ M_r^p ] \leq C^p r^{p/2} p^{2p} \tag{\ref{MART2-3}} \]

for every $p\geq 2$. Furthermore, using inequality \eqref{MART2-3} for $p=2$, we trivially have by Jensen's Inequality

\customlabel{MART2-1}
\[ \E[ M_r^{1/2} ] \leq \E[ M_r^2 ]^{1/4} \leq C^{1/2} r^{1/4} ,   \tag{\ref{MART2-1}} \]

\customlabel{MART2-2}
\[ \E[ M_r^{1} ] \leq \E[ M_r^2 ]^{1/2} \leq C r^{1/2} 2^2     \tag{\ref{MART2-2}}  \]

and

\customlabel{MART2-4}
\[ \hspace{2mm}  \E[ M_r^{3/2} ] \leq \E[ M_r^2 ]^{3/4} \leq C^{3/2} r^{3/4} 2^{3}  .   \tag{\ref{MART2-4}}  \]

Hence, we obtain

\[ \E\left[ \exp\left( \frac18 \left(\frac{M_r}{C \sqrt r}\right)^{1/2} \right) \right] = \sum\limits_{p=0}^\infty 8^{-p} \frac{\E[ M_r^{p/2} ]}{p! C^{p/2} r^{p/4}} . \]

We split the sum for different $p$ and use the above inequalities \eqref{MART2-1}, \eqref{MART2-2}, \eqref{MART2-4} and \eqref{MART2-3} to bound the above expression by

\[ 1 + \underbrace{8^{-1} \frac{C^{1/2} r^{1/4}}{C^{1/2} r^{1/4}}}_{\leq 4^{-1}}  +  \underbrace{8^{-2} \frac{C r^{1/2} 2^2}{C r^{1/2}}}_{=4^{-2}} +  \underbrace{8^{-3} \frac{C^{3/2} r^{3/4} 2^3}{C^{3/2} r^{3/4}}}_{= 4^{-3}}   +  \sum\limits_{p=4}^\infty \underbrace{8^{-p} \frac{(p/2)^p}{p!}}_{\leq 4^{-p} \frac{p^p}{p!}}  \leq 1 + \sum\limits_{p=1}^\infty 4^{-p} \frac{p^p}{p!} .  \]

Using Stirling's Formula for $p\geq 1$

\[ 3 p^p e^{-p} \leq e^{\frac1{12p+1}} \sqrt{2\pi p} p^p e^{-p} \leq p! \]

and the above calculation we finally obtain

\[ \E\left[ \exp\left( \frac18 \left(\frac{M_r}{C \sqrt r}\right)^{1/2} \right) \right] \leq 1 + \frac13 \sum\limits_{p=1}^\infty 4^{-p} e^{p} \leq 2  .  \]

\qed

\end{prof}

\begin{lem}
\xlabel{LEM-EULER}

Let $0 < \eps < \frac16$, $(b_q)_{q\in\N}$ be a sequence of functions $b_q \colon [0,1] \times H \longrightarrow H$ each satisfying Assumption \ref{ASS} then there exists a measurable set $A_{\eps,b} := A_{\eps,(b_q)_{q\in\N}} \subseteq \Omega$, an absolute constant $C \in \R$ and $N_\eps \in \N$ such that for all $x_0 \in Q$, all $n\in\N$ with $n \geq N_\eps$, all $r\in\N$ with $r \leq 2^{n/4}$, $k \in \{ 0 , ..., 2^n - r - 1 \}$ and for every $\eta > 0$ we have

\[ \P \left[  \bbm 1_{A_{\eps,b}^c} \sum\limits_{q=1}^{r}  |\varphi_{n,k+q}(b_q ; x_{q-1}, x_q)|_H   >   \eta C \left( 2^{-n} \sqrt r |x_0|_H + \sqrt r 2^{-2^{n}} \right)   +  C 2^{-n} \sum\limits_{q=0}^{r-1} |x_q|_H  \right]   \leq  4 e^{-\eta^{1/2}} ,  \]

where $x_{q+1} := x_q + \varphi_{n,k+q}(b_q ; x_q)$ for $q\in\{0,...,r-1\}$ is the Euler approximation.

\end{lem}

\begin{prof}

Let $0 < \eps < \frac16$, $n\in\N$ and $b_q\colon [0,1] \times H \longrightarrow H$ be as in the assertion. Using Corollary \ref{COR-SIGMA-RHO} there exists $C_\eps\in\R$ and $A_{\eps,b_q,n,k+q} \in \ca G_{(k+1)2^{-n}}$ with $\P[A_{\eps,b_q,n,k+q}] \leq 2^{-n/24}\eps$ such that for all $x\in 2Q$ we have

\customlabel{THM-EULER-SIGMA}
\[  |\varphi_{n,k+q}(b_q ; x)|_H  \leq  C_\eps n^{\frac12+\frac1\gamma} 2^{-n/2} \left(  |x|_H  +  2^{-2^n} \right)  .  \tag{\ref{THM-EULER-SIGMA}}  \]

on $A_{\eps,b_q,n,k+q}^c$.
Note that $x$ is allowed to be a random variable and we have used that $|\cdot|_\infty \leq |\cdot|_H$.
We now set

\[ N_\eps :=  \min \left\{ n\in\N \setminus \{0\} | C_\eps n^{\frac12+\frac1\gamma} \leq 2^{n/4} \right\}  .  \]

Let, as in the assertion, be $n\in\N$ with $n \geq N_\eps$, $r\leq 2^{n/4}$, $k \in \{ 0 , ..., 2^n - r - 1 \}$ and $x_0 \in Q$.
Additionally, let $x_{q+1} := x_q + \varphi_{n,k+q}(b_q ; x_q)$ be the Euler approximation defined for $q \in \{0, ..., r-1\}$.
We write $x_q = (x_q^{(i)})_{i\in\N}$ for the components of $x_q$ and for $q\in\{1,...,r\}$ we calculate

\[ | x_q^{(i)} |  \leq  | x_{q-1}^{(i)} |  + \left| \int\limits_{(k+q)2^{-n}}^{(k+q+1)2^{-n}} \!\!\!\! b_q^{(i)}(s , Z_s^A+x_q)  -  b^{(i)}_q(s , Z_s^A) \d s \right|   \leq   | x_{q-1}^{(i)} |  +  2 \| b_q^{(i)} \|_\infty 2^{-n} . \]

Via induction on $q$ we deduce

\[ | x_q^{(i)} | \leq | x_0^{(i)} | + 2q 2^{-n} \| b_q^{(i)} \|_\infty \leq 
 |x_0^{(i)}|  + \underbrace{2r 2^{-n}}_{\leq 1} \| b_q^{(i)} \|_\infty  . \]

and since both $x_q\in Q$ and by Assumption \ref{ASS} $b$ takes values in $Q$ we conclude that $x_q \in 2Q$ for all $q \in \{1,...,r\}$.
Note that $x_q$ is $\ca G_{(k+q)2^{-n}}$-measurable.
Due to the fact that Inequality \eqref{THM-EULER-SIGMA} only holds on $A_{\eps,b_q,n,k}^c \subseteq \Omega$ we modify $x_q$ in the following way
\begin{align*}
 \h x_0     &:= x_0 , \\
 \h x_{q+1} &:= \h x_q + \bbm 1_{A_{\eps,b_q,n,k+q}^c} \varphi_{n,k+q}(\h x_q) .
\end{align*}

Observe that we lose the property that $x_{q+1} - x_q = \varphi_{n,k+q}(b_{n,k,q} ; x_q)$, but we still have $\h x_q \in 2Q$ and

\customlabel{THM-EULER-MODIFIED}
\[  |\h x_{q+1} - \h x_q|_H \leq |\varphi_{n,k+q}(b_q ; \h x_q)|_H  .  \tag{\ref{THM-EULER-MODIFIED}}  \]

Furthermore, the modified Euler approximation $\h x_q$ is still $\ca G_{(k+q)2^{-n}}$-measurable. We set

\[ A_{\eps,b} := A_{\eps,(b_q)_{q\in\N}} := \bigcup\limits_{n\in\N} \bigcup\limits_{k=0}^{2^n-1} \bigcup\limits_{q\in\N} A_{\eps,b_q,n,k}  \]

in a similar way as in Corollary \ref{COR-SIGMA-RHO}. We obviously have $\P[A_{\eps,b}] \leq \eps$ and for the modified Euler approximation we obtain for every $q \in \{0, ..., r-1\}$

\[ |\h x_{q+1}|_H = |\h x_q + \bbm 1_{A_{\eps,b_q,n,k+q}^c} \varphi_{n,k+q}(b_{n,k,q} ; \h x_q)|_H   \leq  |\h x_q|_H  +  \bbm 1_{A_{\eps,b_q,n,k+q}^c} |\varphi_{n,k+q}(\h x_q)|_H  \]

and using Inequality \eqref{THM-EULER-SIGMA} for $x$ replaced by $\h x_q$ and $C_\eps n^{\frac12+\frac1\gamma} \leq 2^{n/4}$ this is bounded from above by

\[ |\h x_q|_H + C_\eps n^{\frac12+\frac1\gamma} 2^{-n/2} \left( |\h x_q|_H  +  2^{-2^n} \right)  \leq  (1 + 2^{-n/4}) |\h x_q|_H  +  2^{-n/4} 2^{-2^n}  . \]

By induction over $q \in \{ 0 , ..., r\}$ we have

\[  |\h x_q|_H   \leq  (1 + 2^{-n/4})^q |\h x_0|_H + \sum\limits_{\ell=0}^{q-1} (1 + 2^{-n/4})^\ell 2^{-n/4} 2^{-2^n} .  \]

Using that $q \leq r \leq 2^{n/4}$ this can be further estimated by
   
\[  (1 + 2^{-n/4})^r |x_0|_H + r (1 + 2^{-n/4})^{r} 2^{-n/4} 2^{-2^n}  \leq  \underbrace{(1 + 2^{-n/4})^{2^{n/4}}}_{\leq e} |x_0|_H + \underbrace{r 2^{-n/4}}_{\leq 1} \underbrace{(1 + 2^{-n/4})^{2^{n/4}}}_{\leq e} 2^{-2^n} .  \]

In conclusion we obtain

\customlabel{THM-EULER-x_q}
\[ |\h x_q|_H   \leq  e \left( |x_0|_H   +  2^{-2^n} \right)  .  \tag{\ref{THM-EULER-x_q}}  \]

for all $q \in \{ 0 , ..., r\}$.\\

For the next step we define

\[ \hspace{-53mm}    Y_q := |\varphi_{n,k+q}(b_q ; \h x_{q-1}, \h x_q)|_H  ,  \]
\[                   Z_q := \E [ Y_q  | \ca G_{(k+q)2^{-n}} ] = \E [ |\varphi_{n,k+q}(b_q ; \h x_{q-1}, \h x_q)|_H  | \ca G_{(k+q)2^{-n}} ] , \]
\[ \hspace{-78.5mm}  X_q := Y_q - Z_q ,  \]

and

\[  \hspace{-81mm}  M_\tau := \sum\limits_{q=1}^{r \wedge \tau} X_q .  \]

with $\tau\in\N$. Note that $M_\tau$ is a $\ca G_{(k+\tau+1)2^{-n}}$-Martingale with $M_0 = 0$. Furthermore, for every $p\in\N$ we have the following bound of the increments of $M$

\[ \E[ |X_q|^p ] \leq 2^{p-1} \E [ |Y_q|^p + |Z_q|^p ]  \leq  2^p \E[ |\varphi_{n,k+q}(b_q ; \h x_{q-1}, \h x_q)|_H^p  ] .  \]

Using Corollary \cite[Corollary 3.2]{Wre16} and inequality \eqref{THM-EULER-MODIFIED} this is bounded by

\[ 2^p 3 \beta_A^{p/2} p^{p/2} 2^{-pn/2} \E [ |\h x_{q-1} - \h x_q|_H^p ]  \leq  2^p 3 \beta_A^{p/2} p^{p/2} 2^{-pn/2} \E [ |\varphi_{n,k+q-1}(b_{q-1} ; \h x_{q-1})|_H^p ] . \]

Using Corollary \cite[Corollary 3.2]{Wre16} again this is bounded by

\[ 2^p 9 \beta_A^p p^p 2^{-pn} \E [ |\h x_{q-1} |_H^p ] \leq 18^p \beta_A^p p^p 2^{-pn} \E [ |\h x_{q-1} |_H^p ] . \]

Applying inequality \eqref{THM-EULER-x_q} yields

\[ \E [ |X_q|_H^p ] \leq 18^p \beta_A^p p^p 2^{-pn} e^{p} \left( |x_0|_H   +  2^{-2^{n}}  \right)^p .  \]

Note that $x_0$ is deterministic. Using this bound we invoke Lemma \ref{LEMMA-MARTINGALE-2} with

\[ C := 18 \beta_A 2^{-n} \left( |x_0|_H   +  2^{-2^n} \right) \]

and hence we obtain the following bound for the Martingale $(M_\tau)_{\tau\in\N}$

\customlabel{THM-EULER-M}
\[  \E \left[ \exp \left( \frac18 \left( \frac{r^{-1/2} 2^{n} M_r}{18 \beta_A \left(|x_0|_H  +  2^{-2^n}\right)} \right)^{1/2} \right) \right] \leq 2 .  \tag{\ref{THM-EULER-M}} \]

In a similar way as $(X_q, Y_q, Z_q, M_\tau)$ we define

\[ V_q := \E [ Z_q | \ca G_{(k+q-1)2^{-n}} ] , \]
\[ \hspace{-18mm}  W_q := Z_q - V_q , \]

and

\[ \hspace{-19mm}  M'_\tau := \sum\limits_{\tau=1}^{r \wedge \tau} W_q . \]

Observe that $M'_\tau$ is a $\ca G_{(k+\tau)2^{-n}}$-Martingale and in a completely analogous way as above we obtain

\customlabel{THM-EULER-M'}
\[ \E \left[ \exp \left( \frac18 \left( \frac{r^{-1/2} 2^{n} M'_r}{18 \left(|x_0|_H  +  2^{-2^n}\right)} \right)^{1/2} \right) \right] \leq 2 .  \tag{\ref{THM-EULER-M'}}  \]

Let us now consider the term $V_q$

\[ V_q = \E [ Z_q | \ca G_{(k+q-1)2^{-n}} ]  =  \E [ \E [ |\varphi_{n,k+q}(b_q ; \h x_{q-1},  \h x_{q})|_H   |   \ca G_{(k+q)2^{-n}} ]  |   \ca G_{(k+q-1)2^{-n}} ]   \]

Using Corollary \cite[Corollary 3.2]{Wre16} for $p=1$ and inequality \eqref{THM-EULER-MODIFIED} this is bounded by

\[  3 \beta_A^{1/2} 2^{-n/2} \E [ |\h x_{q-1} -  \h x_{q}|_H  |   \ca G_{(k+q-1)2^{-n}} ]   \leq  3 \beta_A^{1/2} 2^{-n/2} \E [ |\varphi_{n,k+q-1}(b_{q-1} ; \h x_{q-1})|_H |   \ca G_{(k+q-1)2^{-n}} ] . \]

Invoking Corollary \cite[Corollary 3.2]{Wre16} again this can be further bounded from above by

\[  9 \beta_A 2^{-n} \E [ |\h x_{q-1}|_H |  \ca G_{(k+q-1)2^{-n}} ]  =  9 \beta_A 2^{-n} |\h x_{q-1}|_H  .  \]

This leads us to

\customlabel{THM-EULER-V}
\[ \sum\limits_{q=1}^r  V_q  \leq 9 \beta_A 2^{-n} \sum\limits_{q=0}^{r-1} |\h x_{q}|_H   .   \tag{\ref{THM-EULER-V}} \]

For notational ease we set $C' := 18 \beta_A$. Finally, starting from the left-hand side of the assertion and using $Y_q = X_q + W_q + V_q$ we get for every $\eta > 0$

\[ \hspace{4mm}  \P \left[ \bbm 1_{A_{\eps,b}^c} \sum\limits_{q=1}^r \right.   |\varphi_{n,k+q}(b_q ; x_{q-1}, x_q)|_H \left. > \eta C' \left( 2^{-n} \sqrt r |x_0|_H + \sqrt r 2^{-2^{n}} \right)  +  C' 2^{-n} \sum\limits_{q=0}^{r-1} |x_q|_H   \right]  \]

\[ \hspace{-0.5mm}  \leq \P \left[ \sum\limits_{q=1}^r \right.  \bbm 1_{A_{\eps,b}^c} \underbrace{|\varphi_{n,k+q}(b_q ; \h x_{q-1}, \h x_q)|_H}_{=Y_q = X_q + W_q + V_q} \left. > \eta C' \left( 2^{-n} \sqrt r |x_0|_H + \sqrt r 2^{-2^{n}} \right)  +  C' 2^{-n} \sum\limits_{q=0}^{r-1} |x_q|_H   \right]  \]

\[ \hspace{-14mm}  \leq \underbrace{\P \left[ \sum\limits_{q=1}^r V_q  >  C' 2^{-n} \sum\limits_{q=0}^{r-1} |x_q|_H  \right] }_{=0 \text{ by } \eqref{THM-EULER-V}}   +   \P \left[ \sum\limits_{q=1}^r X_q+W_q > \eta C' \sqrt r \left( 2^{-n}  |x_0|_H +  2^{-2^{n}} \right) \right]   \]

\[ \leq   \P \left[ \vphantom{\sum\limits_{q=1}^r} \right. \underbrace{\sum\limits_{q=1}^r X_q}_{=M_r} > C' \eta \sqrt r \left( 2^{-n} |x_0|_H + 2^{-2^{n}} \right) \left. \vphantom{\sum\limits_{q=1}^r} \right]   +  \P \left[ \vphantom{\sum\limits_{q=1}^r} \right. \underbrace{\sum\limits_{q=1}^r W_q}_{=M'_r} > C' \eta \sqrt r \left( 2^{-n} |x_0|_H + 2^{-2^{n}} \right) \left. \vphantom{\sum\limits_{q=1}^r} \right]  \]

\[ = \P \left[  \frac{r^{-1/2} 2^{n}}{C' \left(|x_0|_H + 2^{-2^{n}}\right)} M_r > \eta \right]   +   \P \left[  \frac{r^{-1/2} 2^{n}}{C' \left(|x_0|_H + 2^{-2^{n}} \right)} M'_r > \eta \right] \]

By applying the increasing function $x\mapsto \exp( x^{1/2})$ to both sides and using Chebyshev's Inequality this can be bounded from above by

\[ \exp(- \eta^{1/2} ) \left( \E \left[ \exp\left(  \frac{r^{-1/2} 2^{n}}{C' \left(|x_0|_H + 2^{-2^{n}}\right)} M_r  \right)^{1/2}    +  \exp \left( \frac{r^{-1/2} 2^{n}}{C' \left( |x_0|_H + 2^{-2^{n}} \right)} M'_r \right)^{1/2} \right] \right) . \]

Using inequality \eqref{THM-EULER-M} and \eqref{THM-EULER-M'} we can conclude that

\[ \P \left[ \bbm 1_{A_{\eps,b}^c} \sum\limits_{q=1}^r  |\varphi_{n,k+q}(b_q ; x_{q-1}, x_q)|_H > \eta C' \left( 2^{-n} \sqrt r |x_0|_H + \sqrt r 2^{-2^{n}} \right)  +  C' 2^{-n} \sum\limits_{q=0}^{r-1} |x_q|_H   \right] \leq 4 e^{- \eta^{1/2} } ,  \]

which completes the proof.

\qed

\end{prof}

\begin{thm}
\label{THM-EULER}

For every $0 < \eps < \frac1{40}$ there exist $C_\eps \in \R$, $\Omega_{\eps,b} \subseteq \Omega$ with $\P[\Omega_{\eps,b}^c] \leq \eps$ and $N_\eps\in\N$ such that for all sequences $(b_q)_{q\in\N}$ of functions $b_q \colon [0,1] \times H \longrightarrow H$ with $b_q$ fulfilling Assumption \ref{ASS} for all $q\in\N$, all $n\in\N$ with $n\geq N_\eps$, $k\in  \{0, ..., 2^n - r - 1\}$ and for all $y_0, ..., y_r \in Q$ we have

\[ \sum\limits_{q=1}^r |\varphi_{n,k+q}(b_q ; y_{q-1}, y_q)|_H  \leq C_\eps \left[ 2^{-n} \max \left(r, n^{2+\frac2\gamma} \sqrt r\right) |y_0|_H  + 2^{-n/24} \sum\limits_{q=0}^{r-1} |\gamma_{n,k,q}|_H  + r 2^{-2^{\theta n}} \right] , \]

on $\Omega_{\eps,b}$ for $1 \leq r \leq 2^{n/24}$, where $\gamma_{n,k,q} := y_{q+1} - y_q  - \varphi_{n,k+q}(b_q ; y_q)$ for $q\in\{0,...,r-1\}$ is the error between $y_q$ and the Euler approximation .

\end{thm}

\begin{prof}

\textbf{Step 1:}

Let $0 < \eps < \frac1{40}$ and $C_{\eps/2}$ the constant from Corollary \ref{COR-SIGMA-RHO}. Similar to the proof of Lemma \ref{LEM-EULER} we set

\customlabel{THM-EULER-DE-N}
\[ N_\eps :=  \min \left\{ n\in\N\setminus\{0\} | C_{\eps/2} n^{\frac12+\frac1\gamma} \leq 2^{13n/24} \right\}  .  \tag{\ref{THM-EULER-DE-N}}  \]

For the sake of readability we write $b = (b_q)_{q\in\N}$. By Lemma \ref{LEM-EULER} there is $A_{\eps/2,b} \subseteq \Omega$ with $\P[A_{\eps/2,b}] \leq \frac\eps2$ and a constant $C\in\R$ such that for $x_{q+1} := x_q + \varphi_{n,k+q}(b_q ; x_q)$ and $x_0\in Q$ we have

\customlabel{THM-EULER-LEMMA}
\[   \P \left[ \vphantom{\sum\limits_{q=1}^r} \right. \underbrace{\bbm 1_{A_{\eps/2,b}^c} \sum\limits_{q=1}^{r}  |\varphi_{n,k+q}(b_q ; x_{q-1}, x_q)|_H   >   \eta C \sqrt r \left( 2^{-n}  |x_0|_H  + 2^{-2^{n}} \right)  +  C 2^{-n} \sum\limits_{q=0}^{r-1} |x_q|_H}_{=: B_{\eps/2,b,n,r,k} }  \left. \vphantom{\sum\limits_{q=1}^r} \right]   \leq  4 e^{-\eta^{1/2}}   \tag{\ref{THM-EULER-LEMMA}} \]

for all $\eta > 0$. In order to obtain an almost sure bound we define

\[ B_{\eps/2,b} := \bigcup\limits_{n=N_\eps}^\infty  \bigcup\limits_{r=0}^{2^{n/24}} \bigcup\limits_{k=0}^{2^n-r-1} \bigcup\limits_{s=0}^{2^{2n}} \bigcup\limits_{x_0 \in Q_s \cap 2^{-(s+n)} \Z^\N} B_{\eps/2,b,n,r,k} . \]

Setting

\[ \tilde \eta_\eps := \log \frac{40}\eps \geq 1  \]

and applying Lemma \ref{LEM-EULER} in the form of inequality \eqref{THM-EULER-LEMMA} with $\eta := (1 + 2(3n)^{1+\frac1\gamma})^2 \tilde\eta_\eps^2$ yields

\[ \hspace{-28mm}  \P \left[  B_{\eps/2,b} \right]   \leq  4 \sum\limits_{n=N_\eps}^\infty  \sum\limits_{r=0}^{2^{n/24}} \sum\limits_{k=0}^{2^n-r-1} \sum\limits_{s=0}^{2^{2n}} \sum\limits_{x_0 \in Q_s \cap 2^{-(s+n)} \Z^\N} e^{-\eta^{1/2}}  \]

\[ \hspace{10mm}  \leq  4 \sum\limits_{n=N_\eps}^\infty  2^{n/24} 2^n \sum\limits_{s=0}^{2^{2n}} \#(Q_s \cap 2^{-(s+n)} \Z^\N) e^{-2(3n)^{1+\frac1\gamma}} e^{-\tilde\eta_\eps}  .   \]

Using Theorem \ref{KOLTIK} this is smaller than

\[  4 e^{-\tilde\eta_\eps}  \sum\limits_{n=N_\eps}^\infty  2^{2n} \sum\limits_{s=0}^{2^{2n}} (2\cdot2^n + 1)^{\op{ed}(Q_s)_{s+n}} e^{-2(3n)^{1+\frac1\gamma}} \]

and by invoking Lemma \ref{EFFDIM} this can be again bounded by

\[  \hspace{3mm}  4 e^{-\tilde \eta_\eps}  \sum\limits_{n=N_\eps}^\infty  2^{2n} \sum\limits_{s=0}^{2^{2n}} (2\cdot2^n + 1)^{\ln (s+n+1)^{1/\gamma}} e^{-2(3n)^{1+\frac1\gamma}}  \]

\[  \leq  4 e^{-\tilde \eta_\eps}  \sum\limits_{n=N_\eps}^\infty  2^{4n} (2\cdot 2^n + 1)^{\ln (1 + 2^{2n} + n)^{1/\gamma}} e^{-2(3n)^{1+\frac1\gamma}}   \leq   4 e^{-\tilde \eta_\eps} \sum\limits_{n=N_\eps}^\infty  2^{4n} (3^n)^{ (3n)^{1/\gamma} } e^{-2(3n)^{1+\frac1\gamma}}   \]

\[  \hspace{10mm}  =  4 e^{-\tilde \eta_\eps} \sum\limits_{n=N_\eps}^\infty  2^{4n} \underbrace{3^{(3 n)^{\frac1\gamma}} e^{-(3n)^{1+\frac1\gamma}}}_{\leq1} e^{-(3n)^{1+\frac1\gamma}}   \leq  4 e^{-\tilde \eta_\eps} \underbrace{\sum\limits_{n=N_\eps}^\infty  2^{4n} e^{-3n}}_{\leq 5}   \leq 20 e^{-\tilde \eta_\eps} = \frac\eps2 .  \]

Henceforth, $\P[B_{\eps/2,b}] \leq \frac\eps2$. We set $\Omega_{\eps,b} := A_{\eps/2,b}^c \cap B_{\eps/2,b}^c$. Note that $\P[\Omega_{\eps,b}^c] \leq \eps$.

\medskip

In conclusion there exists $C_\eps\in\R$ such that for all $n\geq N_\eps$, $r \leq 2^{n/24}$, $k \in \{0,...,2^n-r-1\}$ and $x_0 \in Q_s \cap 2^{-(s+n)} \Z^\N$ with $s \in \{ 0, ..., 2^{2n} \}$

\customlabel{THM-EULER-STEP-1}
\[  \sum\limits_{q=1}^r |\varphi_{n,k+q}(b_q ; x_{q-1}, x_q)|_H   \leq   C_\eps \left[ n^{2+\frac2\gamma} 2^{-n} \sqrt r |x_0|_H   +  2^{-n} \sum\limits_{q=0}^{r-1} |x_q|_H   +  r 2^{-2^{\theta n}} \right]  \tag{\ref{THM-EULER-STEP-1}}  \]

holds on $\Omega_{\eps,b}$ with $x_q := x_q + \varphi_{n,k+q}(x_q)$.
Recall that $\theta := \frac23 \frac{\gamma}{\gamma+2}$ and note that we have $\theta \leq \frac23$.\\

\textbf{Step 2:}

Let $n$, $k$, $r\in\N$ and $y_0$, ..., $y_r \in Q$ be as in the statement of this theorem. From now on fix an $\omega\in\Omega_{\eps,b}$. Let $s$ be the largest integer in $\{ 0, ..., 2^{2n} \}$ such that

\[ |y_0|_H \leq 2^{-s} \]

holds. This implies that $y_0 \in Q_s$. Since $s$ is maximal with the above property we have

\[ 2^{-(s+1)} < |y_0|_H   \qquad \text{or} \qquad  |y_0|_H  \leq 2^{-s} = 2^{-2^{2n}} \]

and hence

\customlabel{THM-EULER-2S}
\[ 2^{-s} \leq \max( 2|y_0|_H, 2^{-2^{2n}} ) \leq 2|y_0|_H  +  2^{-2^{2n}} .  \tag{\ref{THM-EULER-2S}} \]

Since $y_0\in Q_s$ we can construct $z_0 \in Q_s \cap 2^{-(s+n)} \Z^\N$, which is close to $y_0$, in the following way: Set $d:=\ln(2s+2n)^{1/\gamma}$. For the components $i  < d$ we choose $z_0$ so that

\customlabel{THM-EULER-FIRST-COMP}
\[ |y_0^{(i)} - z_0^{(i)}| \leq 2^{-s-n}  \tag{\ref{THM-EULER-FIRST-COMP}}  \]

and $z_0^{(i)} := 0$ for $i \geq d$. The distance between $y_0$ and $z_0$ can now be estimated by

\[ |y_0 - z_0|_H^2 = \sum\limits_{0\leq i<d} |y_0^{(i)} - z_0^{(i)}|^2  +  \sum\limits_{d \leq i<\infty}  |y_0^{(i)}|^2  .  \]

Using the above inequality \eqref{THM-EULER-FIRST-COMP} and the fact that $y_0\in Q$ this can be estimated by

\[ d 2^{-2s-2n}  + \sum\limits_{i=d}^{\infty}  4 \exp\left(-2 e^{i^\gamma} \right) \leq d 2^{-2s-2n} + 4 \underbrace{\exp\left(-e^{d^\gamma} \right)}_{=e^{-2s-2n}}  \underbrace{\sum\limits_{i=0}^\infty  \exp\left(-e^{i^\gamma} \right)}_{=: C_\gamma^2 < \infty}  ,  \]

where we have used $\exp(-2e^{i^\gamma}) \leq \exp(-e^{d^\gamma}) \exp(-e^{i^\gamma})$ in the last step. Therefore, we get

\[ |y_0 - z_0|_H \leq 2C_\gamma \sqrt{\ln(2s+2n)^{1/\gamma}} 2^{-s-n} \]

and hence by inequality \eqref{THM-EULER-2S} we obtain

\[ |y_0 - z_0|_H \leq 2C_\gamma \sqrt{\ln(2s+2n)^{1/\gamma}} \left( 2^{1-n} |y_0|_H  +  2^{-n} 2^{-2^{2n}} \right) \]

\[ \leq  4C_\gamma \sqrt{\ln(2^{2n+1}+2n)^{1/\gamma}}  \left( 2^{-n} |y_0|_H  +  2^{-n} 2^{-2^{2n}} \right)   \leq   4C_\gamma \sqrt{\ln(2^{4n})^{1/\gamma}} \left( 2^{-n} |y_0|_H  +  2^{-n} 2^{-2^{2n}} \right)  \]

\[ \hspace{7.5mm}  =  4C_\gamma \sqrt{ \frac{(\log_2(2^{4n}))^{1/\gamma}}{\ln(2)^{1/\gamma}} } \left(  2^{-n} |y_0|_H  +  2^{-n} 2^{-2^n} \right)    =    4C_\gamma \sqrt{ \frac{(4n)^{1/\gamma}}{\ln(2)^{1/\gamma}} } \left( 2^{-n} |y_0|_H  +  2^{-n} 2^{-2^n} \right) . \]

In conclusion we have

\customlabel{THM-EULER-Y_0-Z_0}
\[ |y_0 - z_0|_H  \leq  \tilde C_\gamma \left( n^{\frac1{2\gamma}} 2^{-n} |y_0|_H  +  2^{-2^n} \right)  .  \tag{\ref{THM-EULER-Y_0-Z_0}}  \]

We define $z_1, ..., z_r$ recursively by

\[ z_{q+1} := z_q + \varphi_{n,k+q}(b_q ; z_q) .  \]

Note that $z_0, ..., z_q$ are deterministic since we have fixed $\omega$. Using the definition of $z_q$ we have

\[ |z_{q+1}|_H \leq |z_q|_H  +  |\varphi_{n,k+q}(b_q ; z_q)|_H .  \]

Recall that $\omega\in\Omega_{\eps,b} \subseteq A_{\eps/2,b}^c$ and hence we can invoke the conclusion of Corollary \ref{COR-SIGMA-RHO}, so that the above expression is bounded from above by

\[ |z_q|_H  +  C_{\eps/2} n^{\frac12+\frac1\gamma} 2^{-n/2} ( |z_q|_H + 2^{-2^n})  \leq  (1 + 2^{-n/24}) |z_q|_H  +  2^{-n/24} 2^{-2^n} ,  \]

where we have used Definition \eqref{THM-EULER-DE-N} to conclude that $C_{\eps/2} n^{\frac12+\frac1\gamma} 2^{-n/2} \leq 2^{n/24}$.
By induction on $q \in \{ 1, ..., r-1 \}$ and using $r \leq 2^{n/24}$ we obtain

\[ |z_q|_H \leq (1 + 2^{-n/24})^q  |z_0|_H + \sum\limits_{\ell=0}^{q-1} (1 + 2^{-n/24})^\ell 2^{-n/24} 2^{-2^{n}}  \]

\[ \leq \underbrace{(1 + 2^{-n/24})^r}_{\leq e} |z_0|_H  +  \underbrace{(1 + 2^{-n/24})^r}_{\leq e} \underbrace{r 2^{-n/24}}_{\leq 1} 2^{-2^{n}}  \leq  e \left( |z_0|_H + 2^{-2^{n}} \right) .  \]

Since $z_0$, ..., $z_r$ is by definition an Euler approximation and $z_0 \in Q_s \cap 2^{-(s+n)} \Z^\N$ the conclusion of Step 1 (inequality \eqref{THM-EULER-STEP-1}) with $x_q$ replaced by $z_q$ holds and we obtain that

\[ \hspace{-23mm}  \sum\limits_{q=1}^{r}  |\varphi_{n,k+q}(z_{q-1}, z_q)|_H   \leq   C_\eps \left[ n^{2+\frac2\gamma} 2^{-n} \sqrt r  |z_0|_H  +  r 2^{-2^{\theta n}}  +   2^{-n} \sum\limits_{q=0}^{r-1} |z_q|_H  \right] \]

\[ \hspace{35mm}  \leq C_\eps \left[ n^{2+\frac2\gamma} 2^{-n} \sqrt r  |z_0|_H   +  r 2^{-2^{\theta n}}   +  2^{-n} \sum\limits_{q=0}^{r-1} e ( |z_0|_H + 2^{-2^n} ) \right] \]

\[ \hspace{30mm}  \leq C_\eps \left[ n^{2+\frac2\gamma} 2^{-n} \sqrt r |z_0|_H  +  r 2^{-2^{\theta n}}  +  2^{-n} r e \left( |z_0|_H + 2^{-2^{n}} \right)  \right] \]

\[ \hspace{8mm}  \leq eC_\eps \left[ \max\left(n^{2+\frac2\gamma} \sqrt r, r\right) 2^{-n} |z_0|_H   +  r 2^{-2^{\theta n}} \right]  \]

\[ \hspace{35mm}  \leq  eC_\eps \left[ \max\left(n^{2+\frac2\gamma} \sqrt r, r\right) 2^{-n} ( |y_0|_H + |y_0-z_0|_H )   +  r 2^{-2^{\theta n}} \right] .  \]

Applying inequality \eqref{THM-EULER-Y_0-Z_0} yields that the above expression is bounded from above by

\[ eC_\eps \left[ \max\left(n^{2+\frac2\gamma} \sqrt r, r\right) 2^{-n} \left( |y_0|_H + \tilde C_\gamma \left( \underbrace{n^{\frac1{2\gamma}} 2^{-n}}_{\leq1} |y_0|_H  +  2^{-2^{n}} \right) \right)   +  r 2^{-2^{\theta n}} \right] \]

\[ \hspace{36mm}  \leq  C_{\eps,\gamma} \left[ \max\left(n^{2+\frac2\gamma} \sqrt r, r\right) 2^{-n} |y_0|_H    + r 2^{-2^{\theta n}}  \right]  .  \]

Therefore we obtain

\customlabel{THM-EULER-PHI-Z-Z}
\[ \sum\limits_{q=1}^{r}  |\varphi_{n,k+q}(z_{q-1}, z_q)|_H   \leq  C_{\eps,\gamma} \left[ \max\left(n^{2+\frac2\gamma} \sqrt r, r\right) 2^{-n} |y_0|_H    +  r 2^{-2^{\theta n}} \right] .  \tag{\ref{THM-EULER-PHI-Z-Z}}  \]

\newpage

\textbf{Step 3:}\\

\medskip

\textbf{Claim:}

\customlabel{THM-EULER-PHI-Z-Y}
\[ \sum\limits_{q=1}^r |\varphi_{n,k+q}(z_q, y_q)|_H  \leq   C_\eps' \left[ r 2^{-n} |y_0|_H  +  r 2^{-2^{\theta n}}   +  2^{-n/24} \sum\limits_{q=0}^{r-1} | \gamma_{n,k,q} |_H  \right] .  \tag{\ref{THM-EULER-PHI-Z-Y}}  \]

\textbf{Proof of \eqref{THM-EULER-PHI-Z-Y}}:

\medskip

We set $u_q := z_q - y_q$ for $q \in \{0, ..., r \}$ and bound the increments of $u_q$ in the following way.

\[ |u_{q+1} - u_q|_H = |z_{q+1} - y_{q+1} - z_q + y_q|_H = | \varphi_{n,k+q}(b_q ; z_q) - y_{q+1} + y_q |_H \]

\[ \hspace{8mm}  \leq | \varphi_{n,k+q}(b_q ; z_q) - y_{q+1} + y_q + \gamma_{n,k,q} |_H  +  | \gamma_{n,k,q} |_H \]

\[ \hspace{2mm}  = | \varphi_{n,k+q}(b_q ; z_q) - \varphi_{n,k+q}(b_q ; y_q) |_H  +  | \gamma_{n,k,q} |_H .  \]

\[ \hspace{-22mm} = |\varphi_{n,k+q}(b_q ; z_q, y_q)|_H  +  | \gamma_{n,k,q} |_H \]

We therefore deduce that

\[ |u_{q+1}|_H \leq |u_{q+1} - u_q|_H  +  |u_q|_H \leq | \varphi_{n,k+q}(b_q ; z_q, y_q) |_H  +  | \gamma_{n,k,q} |_H  +  |u_q|_H  .  \]

By the conclusion of Corollary \ref{COR-SIGMA-RHO} and Definition \eqref{THM-EULER-DE-N} this is bounded by

\[ C_{\eps/2} \left( \vphantom{2^{-2^{\theta n}}} \right. \sqrt{n} 2^{-n/6} \underbrace{|z_q-y_q|_H}_{=|u_q|_H}  + 2^{-2^{\theta n}} \left. \vphantom{2^{-2^{\theta n}}} \right)  +   | \gamma_q |_H + |u_q|_H  \leq  ( 1 + 2^{-n/24} ) |u_q|_H  + C_{\eps/2} 2^{-2^{\theta n}}   +  | \gamma_{n,k,q} |_H   .  \]

Induction on $q \in \{0,...,r\}$ yields

\[ |u_q|_H  \leq  C_{\eps/2} ( 1 + 2^{-n/24} )^r \left( |u_0|_H  + r 2^{-2^{\theta n}}   +  \sum\limits_{q=0}^{r-1} | \gamma_{n,k,q} |_H \right)  \]

\[ \hspace{-9mm}  \leq  e C_{\eps/2} \left( |u_0|_H  + r 2^{-2^{\theta n}}   +  \sum\limits_{q=0}^{r-1} | \gamma_{n,k,q} |_H  \right)  .  \]

Using inequality \eqref{THM-EULER-Y_0-Z_0} together with the above calculation yields

\[ |u_q|_H  \leq  e C_{\eps/2}  \left( \tilde C_\gamma n^{\frac1{2\gamma}} 2^{-n} |y_0|_H  +  2r 2^{-2^{\theta n}}   +  \sum\limits_{q=0}^{r-1} | \gamma_{n,k,q} |_H  \right) \]

and hence by combining this estimate with Corollary \ref{COR-SIGMA-RHO} we have

\[ \hspace{-10mm}  |\varphi_{n,k+q}(z_q, y_q)|_H \leq C_{\eps/2} \left( \sqrt n 2^{-n/6} |z_q-y_q|_H  +  2^{-2^{\theta n}} \right)  \leq  C_{\eps/2} \left( 2^{-n/12} |u_q|_H  +  2^{-2^{\theta n}} \right) \]

\[ \hspace{27mm}  \leq eC_{\eps/2}^2 2^{-n/12} \left( \tilde C_\gamma n^{\frac1{2\gamma}} 2^{-n} |y_0|_H  +  2r 2^{-2^{\theta n}}    +  \sum\limits_{q=0}^{r-1} | \gamma_{n,k,q} |_H   \right)  +  C_{\eps/2} 2^{-2^{\theta n}}  .  \]

In conclusion since $r\leq 2^{n/24}$ we obtain

\[ |\varphi_{n,k+q}(z_q, y_q)|_H  \leq   C_\eps' \left[ 2^{-n} |y_0|_H  +  2^{-2^{\theta n}}   +  2^{-n/12} \sum\limits_{q=0}^{r-1} | \gamma_{n,k,q} |_H  \right]    \]

and hence summing over $q=1,...,r$ and using again that $r\leq 2^{n/24}$ complete the proof of Claim \eqref{THM-EULER-PHI-Z-Y}.

\medskip

\textbf{Step 4:}

\medskip

Finally, using the identity $y_{q-1} - y_q = y_{q-1} - z_{q-1} + z_{q-1} - z_q + z_q - y_q$ the left-hand side of the assertion can be bounded as follows

\[ \!\! \sum\limits_{q=1}^r \! |\varphi_{n,k+q}( b_q ; y_{q-1}, y_q )|_H  \leq  \sum\limits_{q=1}^r \! |\varphi_{n,k+q}( b_q ; y_{q-1}, z_{q-1} )|_H  +  |\varphi_{n,k+q}( b_q ; z_{q-1}, z_q )|_H   +   |\varphi_{n,k+q}( b_q ; z_q, y_q )|_H  .  \]

Applying inequalities \eqref{THM-EULER-PHI-Z-Z}, \eqref{THM-EULER-PHI-Z-Y} and \eqref{THM-EULER-PHI-Z-Y} with $z_q$, $y_q$ replaced by $z_{q-1}$, $y_{q-1}$ respectively yields that this is bounded by

\[ C_\eps'' \left[ 2^{-n} \max \left(r, n^{2+\frac2\gamma} \sqrt r\right) |y_0|_H   +  r 2^{-2^{\theta n}}   +  2^{-n/24} \sum\limits_{q=0}^{r-1} | \gamma_{n,k,q} |_H  \right] .  \]

\qed

\end{prof}

\begin{cor}
\xlabel{COR-GLUING}

For every $0 < \eps < \frac1{40}$ there exists $C_\eps \in \R$ such that for every sequence $(b_q)_{q\in\N}$ of Borel measurable functions $b_q \colon [0,1] \times H \longrightarrow H$ satisfying Assumption \ref{ASS} there exists a measurable set $\Omega_{\eps,(b_q)_{q\in\N}} \subseteq \Omega$ with $\P[\Omega_{\eps,(b_q)_{q\in\N}}^c] \leq \eps$ such that for all sufficiently large $n\in\N$, $N \in\N$ with $N \leq 2^n$, $k\in \{0, ..., 2^n - N \}$ and for all $x_q \in Q$ for $q \in \{ 0 ,..., N \}$, we have

\[ \sum\limits_{q=0}^{N-1} |\varphi_{n,k+q}(b_q ; x_{q+1}, x_q)|_H  \leq C_\eps \left[ 2^{-n} \sum\limits_{q=0}^{N} |x_q|_H   +   2^{-3n/4} |x_0|_H  + 2^{-n/24} \sum\limits_{q=0}^{N-1} |\gamma_{n,k,q}|_H  + N 2^{-2^{\theta n}} \right] , \]

on $\Omega_{\eps,b}$, where $\gamma_{n,k,q} := x_{q+1} - x_q  - \varphi_{n,k+q}(b_q ; x_q)$ is the error between $x_{q+1}$ and the Euler approximation for $x_{q+1}$ given $x_q$.

\end{cor}

\begin{prof}

We set $r := \lfloor 2^{n/24} \rfloor$. For the sake of notional ease we set $x_{q'} = 0$ whenever $q'>N$. In order to estimate the left-hand side of the assertion we will use Theorem \ref{THM-EULER}. To this end we split the sum into $s$ pieces of size $r$. Choose $i\in\{ 0, ..., r-1 \}$ such that

\[  \sum\limits_{t=0}^{\lfloor r^{-1} N \rfloor} |x_{i + tr}|_H \leq \frac1r \sum\limits_{q=0}^{r-1} \sum\limits_{t=0}^{\lfloor r^{-1} N \rfloor}  |x_{q +tr}|_H \]

holds. Since we calculate the mean of $\sum\limits_{t=0}^{\lfloor r^{-1} N \rfloor} |x_{q + tr}|_H$ on the right-hand side, it is clear that such an $i$ always exists. Set $s := \lfloor r^{-1} (N - i) \rfloor$ and note that $s\leq \lfloor r^{-1} N \rfloor$. Using this we have

\[ \sum\limits_{t=0}^s |x_{i + tr}|_H  \leq  \frac1r \sum\limits_{q=0}^{r-1} \sum\limits_{t=0}^{\lfloor r^{-1} N \rfloor} |x_{q +tr}|_H .  \]

Hence, we obtain

\customlabel{COR-GLUING-PSI}
\[ \sum\limits_{t=0}^s |x_{i + tr}|_H \leq r^{-1} \sum\limits_{q=0}^{N-1} |x_q|_H .  \tag{\ref{COR-GLUING-PSI}}  \]

Starting with the left-hand side of the assertion we split the sum into three parts. The first part contains the terms $x_q$ for $q=0$ to $q=i$. Since $i \leq r \leq 2^{n/24}$ this can be handled by applying Theorem \ref{THM-EULER} directly. The second part contains $s$ sums of size $r$. Here, Theorem \ref{THM-EULER} is applicable for every term of the outer sum running over $t$. The last part can be handled, in the same way as the first part, by directly applying Theorem \ref{THM-EULER}. This strategy yields

\[  \hspace{-50mm}  \sum\limits_{q=0}^{N-1}  |\varphi_{n,k+q}( b_q ; x_{q+1}, x_q )|_H  =  \sum\limits_{q=0}^{i-1}  |\varphi_{n,k+q}( b_q ; x_{q+1}, x_q )|_H \]
\[ \hspace{10mm}  \qquad\qquad\qquad  +  \sum\limits_{t=0}^{s-1} \sum\limits_{q=0}^{r-1}  |\varphi_{n,k + i + tr + q}( b_q ; x_{q+1 + i + tr}, x_{q + i + tr} )|_H  \]
\[ \hspace{13mm}  \qquad\qquad\qquad  +  \sum\limits_{q=0}^{N-i-rs-1}  |\varphi_{n,k + i + sr +q}( b_q ; x_{q+1 + i + sr}, x_{q + i + sr} )|_H \]

\newpage

\[ \hspace{-25mm} \leq C_\eps \left[ 2^{-n} \max \left(r, n^{2+\frac2\gamma} \sqrt r\right) |x_0|_H  +  2^{-n/24} \sum\limits_{q=0}^{i-1} |\gamma_{n,k,q}|_H  +  r 2^{-2^{\theta n}} \right] \]

\[  \hspace{-24mm}  + C_\eps \sum\limits_{t=0}^{s-1} \left[ 2^{-n} r |x_{i + tr}|_H + 2^{-n/24}  \sum\limits_{q=0}^{r-1} |\gamma_{n,k,i + tr + q}|_H  +  r 2^{-2^{\theta n}} \right]  \]

\[  \hspace{-1mm}  \quad\qquad + C_\eps \left[ 2^{-n} \max \left(r, n^{2+\frac2\gamma} \sqrt r\right) |x_{i + sr}|_H  +  2^{-n/24}  \sum\limits_{q=0}^{N-i-rs-1} |\gamma_{n,k,i + sr + q}|_H  +  r 2^{-2^{\theta n}} \right] .  \]

\[ \hspace{-7mm} \leq  C_\eps \left[ 2^{-n} r |x_0|_H  +  2^{-n} r \sum\limits_{t=0}^s |x_{i + tr}|_H + 2^{-n/24} \sum\limits_{q=0}^{N-1} |\gamma_{n,k,q}|_H   +  (s+2) r 2^{-2^{\theta n}} \right] .  \]

Estimating this further by using inequality \eqref{COR-GLUING-PSI} and $r \leq 2^{n/24}$ yields the following bound

\[ 2 C_\eps \left[ 2^{-3n/4} |x_0|_H  +  2^{-n} \sum\limits_{q=0}^{N-1} |x_q|_H + 2^{-n/24} \sum\limits_{q=0}^{N-1} |\gamma_{n,k,q}|_H   +  N 2^{-2^{\theta n}} \right] ,  \]

which completes the proof.

\qed

\end{prof}

\section{Proof of the main result}

In this section we are going to formulate a $\log$-type Gronwall inequality of the form

\[ f((j+1)2^{-n}) - f(j2^{-n}) \leq C 2^{-n} f(j2^{-n}) \log(1/f(j2^{-n})) \]
\[ \Rightarrow f(j2^{-n}) \leq C f(0)    \]

for $j \in \{0, ..., 2^n\}$.

\medskip

In Lemma \ref{LEM-GRONWALL} we prove this implication in an abstract setting.
Using all our previous considerations we show in Theorem \ref{THM-FINAL} that our function $u$ from Proposition \ref{PRO-GIRSANOV} satisfies such a Gronwall inequality and hence has to coincide with the zero function (Corollary \ref{COR-FINAL}).

\begin{lem}[$\log$-Type Gronwall Inequality]
\xlabel{LEM-GRONWALL}

Let $K>0$, $m\in\N$ ``sufficiently big'' i.e.~$K \leq \ln(2) 2^{m}$ and $0 < \beta_0, ..., \beta_{2^m} < 1$ and assume that

\[ \Delta \beta_j \leq K 2^{-m} \beta_j \log_2(1/\beta_j), \qquad \forall j \in \{ 0, ..., 2^m-1 \} \]

holds, where $\Delta \beta_j := \beta_{j+1} - \beta_j$. Then, we have

\[ \beta_j \leq \exp\left( \log_2(\beta_0) e^{-2K-1} \right), \qquad \forall j \in \{ 0, ..., 2^m \} . \]

\end{lem}

\begin{prof}

For every $j \in \{ 0, ..., 2^m \}$ we define

\[ \gamma_j := \log_2(1/\beta_j) . \]

By assumption we have

\[ \hspace{-34mm}  \gamma_{j+1} = - \log_2(\beta_{j+1}) \geq - \log_2( \beta_j + K 2^{-m} \beta_j \gamma_j) \]

\[ \hspace{12mm}  = - \log_2(\beta_j) - \log_2( 1 + K 2^{-m} \gamma_j)
= \gamma_j - \frac1{\ln 2} \ln( 1 + K 2^{-m} \gamma_j)  .  \]

Using the inequality $\ln(1+x) \leq x$ the above, and hence $\gamma_{j+1}$, is smaller than

\[ \gamma_j \left( 1 - \frac{K}{\ln 2} 2^{-m} \right) . \]

By induction on $j\in \{0, ..., 2^m\}$ we obtain

\[ \gamma_j \geq \gamma_0 \left( 1 - \frac{K}{\ln 2} 2^{-m} \right)^j .  \]

Since $m$ is ``sufficiently big'' the term inside the brackets is in the interval $[0,1]$ so that $\gamma_j$ is bounded from below by

\[ \gamma_0 \left( 1 - \frac{K}{\ln 2} 2^{-m} \right)^{2^m} \geq \gamma_0 e^{- K/\ln(2) - 1} \geq \gamma_0 e^{- 2K - 1} .  \]

Plugging in the definition of $\gamma_j$ implies that

\[ \log_2(1/\beta_j)  \geq  \log_2(1/\beta_0) e^{- 2K - 1} . \]

Isolating $\beta_j$ yields

\[ \beta_j  \leq  \exp\left( \log_2(\beta_0) e^{- 2K - 1} \right) . \]

\qed

\end{prof}

\begin{thm}
\xlabel{THM-FINAL}

Let $0 < \eps < \frac1{40}$ and $f$ be as in Assumption \ref{ASS} then there exist $A_{\eps,f} \subseteq \Omega$, $K=K(\eps)>0$ and $m_0 = m_0(\eps)\in\N$ with $\P[A_{\eps,f}^c ]\leq \eps$ such that for any function $u\in\Phi$ being a solution of equation \eqref{DE-U} for a fixed $\omega\in A_{\eps,f}$, for all integers $m$ with $m \geq m_0$, $j\in\{0,...,2^m-1\}$ and $\beta$ satisfying

\[ 2^{m-2^{\theta m}} \leq \beta \leq 2^{-2^{(\frac\theta2+\frac14)m}} \quad \text{and} \quad |u(j2^{-m})|_H \leq \beta \]

we have

\[ |u((j+1)2^{-m})|_H \leq \beta \left( 1 + K2^{-m} \log_2(1/\beta) \right) .  \]

\end{thm}

\begin{prof}

Let $0 < \eps < \frac1{40}$ and $f$ be as in the assertion. For all $n\in\N$ and $k\in\{0, ..., 2^n-1\}$ we set

\[ b_{n,k}(t, x) := e^{-((k+1)2^{-n}-t)A} f(t,x)  , \qquad \forall t \in [0,1] , \ x \in H .  \]

Note that $b_{n,k}$ fulfills Assumption \ref{ASS} since $|b_{n,k}(t,x)|_H \leq |f(t,x)|_H$. Choose $A_{\eps,f} \subseteq \Omega$ with $\P[A_{\eps,f}] \leq \eps$ such that the conclusions of Corollary \ref{COR-SIGMA-RHO}, Theorem \ref{THM-APPROX} and Corollary \ref{COR-GLUING} hold with the same constant $C_\eps\geq 1$ for all functions $b_{n,k}$ on $A_{\eps,f}^c$. We set

\[ m_0 := \max\left( 3\log_2(584C_\eps'), 24\log_2(72C_\eps) , \frac{1}{2\theta-1} \right) , \]

where $C_\eps'$ will be defined later. Recall that we defined $\theta := \frac23 \frac\gamma{\gamma+2}$. Fix an $\omega\in A_{\eps,f}^c$, $m\geq m_0$, $j$, $u$ and $\beta$ as in the statement and suppose $|u(j2^{-m})|_H \leq \beta$. We set $N := 7 \lfloor \log_2(1/\beta) \rfloor$. Observe that

\customlabel{THM-FINAL-DE-N}
\[ m^{\frac12+\frac1\gamma} 2^{m/2} \leq 7 \cdot 2^{(\frac\theta2+\frac14)m} - 7 \leq N \leq 7 \cdot 2^{\theta m} \leq 7 \cdot 2^{2m/3} ,  \tag{\ref{THM-FINAL-DE-N}}  \]

where we have used that $\frac12 < \theta \leq \frac23$ (due to $\gamma>6$ in Assumption \ref{ASS} and $\theta := \frac23 \frac{\gamma}{\gamma+2}$).
Suppose $u\in\Phi$ satisfies equation \eqref{DE-U} as stated in the assertion. We define for every $n\in\N$ and $t\in[0,1]$

\[ u_n(t) := \sum\limits_{k=0}^{2^n-1} \bbm 1_{[k2^{-n}, (k+1)2^{-n}[}(t) u(k2^{-n}) .  \]

Note that $u_n$ converges pointwise to $u$ on $[0,1[$ and $u_n\in\Phi^*$ by construction and since $u\in\Phi$. Let $\alpha$ be the smallest real number such that

\customlabel{THM-FINAL-DE-ALPHA}
\[ \sum\limits_{k=j2^{n-m}}^{(j+1)2^{n-m}-1} \!\!\!\!\! |u((k+1)2^{-n}) - u(k2^{-n})|_H \leq \alpha 2^{-m} \left[ N + n^{\frac12+\frac1\gamma} 2^{n/2} \right],  \quad \forall n\in\{ m, ..., N\}   \tag{\ref{THM-FINAL-DE-ALPHA}}  \]

holds. I.e.

\[ \alpha := \max\limits_{m \leq n \leq N} \frac{2^m}{N + n^{\frac12+\frac1\gamma} 2^{n/2}}   \sum\limits_{k=j2^{n-m}}^{(j+1)2^{n-m}-1} \!\!\! |u((k+1)2^{-n}) - u(k2^{-n})|_H  .  \]

For $n\geq m$ we define

\customlabel{THM-FINAL-DE-PSI}
\[ \psi_n := \sum\limits_{k=j2^{n-m}}^{(j+1)2^{n-m}-1} \!\! |u(k2^{-n})|_H .    \tag{\ref{THM-FINAL-DE-PSI}}  \]

By splitting the sum in \eqref{THM-FINAL-DE-PSI} in two sums, one where $k$ is even and one where $k$ is odd, we can estimate $\psi_n$ by $\psi_{n-1}$.
To this end let $n\in\{m+1,...,N\}$. We then have

\[  \psi_n = \sum\limits_{\bfrac{k=j2^{n-m}}{2\mid k}}^{(j+1)2^{n-m}-1}   |u(k2^{-n})|_H  +  \!\!\! \sum\limits_{\bfrac{k=j2^{n-m}}{2\nmid k}}^{(j+1)2^{n-m}-1} \!\!\!\! |u(k2^{-n})|_H \]

\[  \leq  \!\!\!\!\!\! \sum\limits_{\bfrac{k=j2^{n-m}}{2\mid k}}^{(j+1)2^{n-m}-1} \!\!\!\!\!\!\!\! |u(k2^{-n})|_H  +  \!\!\!\!\!\!\! \sum\limits_{\bfrac{k=j2^{n-m}}{2\nmid k}}^{(j+1)2^{n-m}-1} \!\!\!\!\!\!\!\! |u(k2^{-n}) - u((k-1)2^{-n})|_H + |u((k-1)2^{-n})|_H + |u((k+1)2^{-n}) - u(k2^{-n})|_H  .  \]

Since $k-1$ is even whenever $k$ is odd, rewriting the term $|u((k-1)2^{-n})|_H$ yields that the above equals

\[ \hspace{-84mm}  \sum\limits_{\bfrac{k=j2^{n-m}}{2\mid k}}^{(j+1)2^{n-m}-1}   |u( k 2^{-n})|_H  + |u( k 2^{-n})|_H  \]

\[  \hspace{10mm}    +   \sum\limits_{\bfrac{k=j2^{n-m}}{2\nmid k}}^{(j+1)2^{n-m}-1}  |u(k2^{-n}) - u((k-1)2^{-n})|_H + |u((k+1)2^{-n}) - u(k2^{-n})|_H  \]

\[   = 2  \!\!\! \sum\limits_{k=j2^{n-m-1}}^{(j+1)2^{n-m-1}-1} \!\!\!\! |u(k 2^{-n+1})|_H   +   \! \sum\limits_{\bfrac{k=j2^{n-m}}{2\nmid k}}^{(j+1)2^{n-m}-1} \!\!\! |u(k2^{-n}) - u((k-1)2^{-n})|_H + |u((k+1)2^{-n}) - u(k2^{-n})|_H  \]

\[  \hspace{-32mm}   = 2 \! \sum\limits_{k=j2^{n-1-m}}^{(j+1)2^{n-1-m}-1} \! |u(k 2^{-(n-1)})|_H  +  \sum\limits_{k=j2^{n-m}}^{(j+1)2^{n-m}-1} \!\! |u((k+1)2^{-n}) - u(k2^{-n})|_H  .  \]

And, henceforth, since $n\in \{m+1,...,N\}$ using inequality \eqref{THM-FINAL-DE-ALPHA} and the definition of $\psi_n$ (equation \eqref{THM-FINAL-DE-PSI}) we have the following bound

\[  \psi_n \leq 2 \psi_{n-1}  +  \alpha 2^{-m} N  +  \alpha 2^{-m} n^{\frac12+\frac1\gamma} 2^{n/2}  .  \]

By induction we deduce

\[  \hspace{-42mm}  \psi_n \leq 2^{n-m} \psi_m  +  \sum\limits_{\ell=m+1}^n \alpha 2^{n-\ell-m} N  +   \sum\limits_{\ell=m+1}^n \alpha 2^{n-\ell-m} \ell^{\frac12+\frac1\gamma} 2^{\ell/2}  \]

\[ \hspace{11mm}  \leq  2^{n-m} |u(j2^{-m})|_H  +  \alpha 2^{n-m} N \sum\limits_{\ell=m+1}^n 2^{-\ell}  +  \alpha 2^{n-m} \sum\limits_{\ell=m+1}^n 2^{-\ell/3} ,  \qquad \forall n\in\{m,...,N\} .  \]

We use $|u(j2^{-m})|_H \leq \beta$ to bound the above by

\[ 2^{n-m} \left[ \beta  + \alpha N \sum\limits_{\ell=m+1}^n 2^{-\ell} + \alpha \sum\limits_{\ell=m+1}^n 2^{-\ell/3} \right] \leq 2^{n-m} \left[ \beta  +  \alpha 2^{-m} N + \alpha 2^{-m/3} \right] .  \]

Furthermore, using inequality \eqref{THM-FINAL-DE-N} i.e.~$N \leq 7\cdot 2^{2m/3}$ we bound the above by

\[ 2^{n-m} \left[ \beta  +  7 \alpha 2^{-m/3}  +  \alpha 2^{-m/3} \right]  \leq 8 \cdot 2^{n-m} \left[ \beta  +  \alpha 2^{-m/3} \right]  .  \]

In conclusion we obtain

\customlabel{THM-FINAL-PSI-ESTIMATE}
\[ \psi_n \leq 8\cdot 2^{n-m} \left( \beta  +  \alpha 2^{-m/3} \right), \qquad\qquad \forall n\in\{m,...,N\} .  \tag{\ref{THM-FINAL-PSI-ESTIMATE}}  \]

Since $u$ solves equation \eqref{DE-U} we have

\[  \left| u((k+1)2^{-n}) - u(k2^{-n}) - \varphi_{n,k}(b_{n,k} ; u(k2^{-n})) \right|_H \]

\[ = \left| u((k+1)2^{-n}) - u(k2^{-n}) - \int\limits_{k2^{-n}}^{{(k+1)2^{-n}}}  b_{n,k}(t,Z_t^A(\omega) + u(k2^{-n})) - b_{n,k}(t,Z_t^A) \d t \right|_H \]

\[ \hspace{-33mm}  \overset{\eqref{DE-U}}= \left| \int\limits_0^{(k+1)2^{-n}} e^{-((k+1)2^{-n}-t)A} ( f(t, Z^A_t(\omega) + u(t)) - f(t, Z^A_t(\omega)) ) \d t \right. \]

\[ \hspace{-25mm}  -   \int\limits_0^{k2^{-n}} e^{-(k2^{-n}-t)A} (f(t, Z^A_t(\omega) + u(t)) - f(t, Z^A_t(\omega))) \d t   \]

\[ \hspace{-20mm}  -   \left. \int\limits_{k2^{-n}}^{{(k+1)2^{-n}}} b_{n,k}(t,Z_t^A(\omega) + u(k2^{-n}))  -  b_{n,k}(t,Z_t^A(\omega)) ) \d t \right|_H   \]

\[ \hspace{-32mm}  = \left| \int\limits_{k2^{-n}}^{(k+1)2^{-n}} e^{-((k+1)2^{-n}-t)A} (f(t, Z^A_t(\omega) + u(t)) - f(t, Z^A_t(\omega))) \d t   \right. \]

\[ \hspace{-28mm} -   \int\limits_{k2^{-n}}^{{(k+1)2^{-n}}} b_{n,k}(t,Z_t^A(\omega) + u(k2^{-n}))  - b_{n,k}(t,Z_t^A(\omega)) ) \d t   \]

\[ \hspace{18mm}  +  \left. \int\limits_0^{k2^{-n}} \left( e^{-((k+1)2^{-n} - t)A} -  e^{-(k2^{-n} - t)A} \right) \cdot ( f(t,Z_t^A(\omega) + u(t))  - f(t,Z_t^A(\omega)) ) \d t   \right|_H . \]

Using the definition of $b_{n,k}$ this can be simplified and bounded by

\[ \hspace{-32mm}  \left| \int\limits_{k2^{-n}}^{(k+1)2^{-n}}  b_{n,k}(t, Z^A_t(\omega) + u(t)) - b_{n,k}(t, Z^A_t(\omega) + u(k2^{-n}))   \d t \right|_H \]

\[ \hspace{5mm}  +  \underbrace{\left| e^{-2^{-n} A} - 1 \right|_{\text{op}}}_{\leq C2^{-n}} \cdot \left| \int\limits_0^{k2^{-n}} e^{-(k2^{-n} - t)A} ( f(t,Z_t^A(\omega) + u(t))  - f(t,Z_t^A(\omega)) ) \d t \right|_H . \]

Since $u_n$ is constant on $[k2^{-n}, (k+1)2^{-n}[$ and using again that $u$ solves equation \eqref{DE-U} we can estimate this by

\[  \left| \int\limits_{k2^{-n}}^{(k+1)2^{-n}} b_{n,k}(t, Z^A_t(\omega) + u(t)) - b_{n,k}(t, Z^A_t(\omega) + u(k2^{-n}))  \d t \right|_H   +  C 2^{-n} |u(k2^{-n})|_H . \]

By invoking Theorem \ref{THM-APPROX} this can be rewritten as

\[ \lim\limits_{\ell\rightarrow\infty} \left| \int\limits_{k2^{-n}}^{(k+1)2^{-n}}  b_{n,k}(t, Z^A_t(\omega) + u_\ell(t)) - b_{n,k}(t, Z^A_t(\omega) + u_n(t))  \d t \right|_H +  C 2^{-n} |u(k2^{-n})|_H \]

\[ \!\!\!\! \!\! \!\!\!\!\!\! \leq  C 2^{-n} |u(k2^{-n})|_H  +  \sum\limits_{\ell=n}^\infty \left| \int\limits_{k2^{-n}}^{(k+1)2^{-n}} \! b_{n,k}(t, Z^A_t(\omega) + u_{\ell+1}(t)) - b_{n,k}(t, Z^A_t(\omega) + u_\ell(t)) \d t \right|_H .  \]

\[ \!\!\!\! =  C 2^{-n} |u(k2^{-n})|_H  +  \sum\limits_{\ell=n}^\infty \! \sum\limits_{r=k2^{\ell-n}}^{(k+1)2^{\ell-n}-1} \left| \!\int\limits_{2r2^{-\ell-1}}^{(2r+2)2^{-\ell-1}}\!\!\! b_{n,k}(t, Z^A_t(\omega) + u_{\ell+1}(t)) - b_{n,k}(t, Z^A_t(\omega) + u_\ell(t)) \d t \right|_H \]

\[  \!\!\!\! =  C 2^{-n} |u(k2^{-n})|_H  +  \sum\limits_{\ell=n}^\infty \!\! \sum\limits_{r=k2^{\ell-n}}^{(k+1)2^{\ell-n}-1} \left| \!\int\limits_{2r2^{-\ell-1}}^{(2r+1)2^{-\ell-1}}\!\!\!\!\!\!\!\!\! \right. \underbrace{b_{n,k}(t, Z^A_t(\omega) + u(2r2^{-\ell-1})) - b_{n,k}(t, Z^A_t(\omega) + u(r2^{-\ell}))}_{=0} \!\d t \left. \vphantom{\int\limits_{2r2^{-\ell-1}}^{(2r+1)2^{-\ell-1}}}\right|_H \]

\[ + \sum\limits_{\ell=n}^\infty  \sum\limits_{r=k2^{\ell-n}}^{(k+1)2^{\ell-n}-1}   \left| \ \int\limits_{(2r+1)2^{-\ell-1}}^{(2r+2)2^{-\ell-1}}\!\!\!\!\! b_{n,k}(t, Z^A_t(\omega) + u((2r+1)2^{-\ell-1})) - b_{n,k}(t, Z^A_t(\omega) + u(r2^{-\ell})) \d t \!\right|_H \]

\[  =  C 2^{-n} |u(k2^{-n})|_H  +  \sum\limits_{\ell=n}^\infty \! \sum\limits_{r=k2^{\ell-n}}^{(k+1)2^{\ell-n}-1} \!\!\! |\varphi_{\ell+1,2r+1}\left( b_{n,k} ; u\left((2r+1)2^{-\ell-1}\right), u\left(r2^{-\ell}\right) \right)|_H .  \]

Summing over $k\in\{j 2^{n-m}, ..., (j+1)2^{n-m}-1 \}$ leads us to

\[ \hspace{-49mm}  \sum\limits_{k=j2^{n-m}}^{(j+1)2^{n-m}-1} \!\! |u((k+1)2^{-n}) - u(k2^{-n}) - \varphi_{n,k}(b_{n,k} ; u(k2^{-n}))|_H \]

\[ \leq  \!\!\! \sum\limits_{k=j2^{n-m}}^{(j+1)2^{n-m}-1} \! \left( C 2^{-n} |u(k2^{-n})|_H  + \sum\limits_{\ell=n}^\infty \! \sum\limits_{r=k2^{\ell-n}}^{(k+1)2^{\ell-n}-1} \!\!\!\! |\varphi_{\ell+1,2r+1}\left( b_{n,k} ; u\left((2r+1)2^{-\ell-1}\right), u\left(r2^{-\ell}\right) \right)|_H \right) .  \]

\[ = \!\!\!\!\! \sum\limits_{k=j2^{n-m}}^{(j+1)2^{n-m}-1} \!\!\!\!\!\!\! C 2^{-n} |u(k2^{-n})|_H  +  \sum\limits_{\ell=n}^\infty \! \sum\limits_{k=j2^{n-m}}^{(j+1)2^{n-m}-1} \sum\limits_{r=k2^{\ell-n}}^{(k+1)2^{\ell-n}-1} \!\!\!\!\!\! |\varphi_{\ell+1,2r+1}\left(b_{n,k} ; u\left((2r+1)2^{-\ell-1}\right), u\left(r2^{-\ell}\right)\right)|_H  \]

\[ = \!\!\! \sum\limits_{k=j2^{n-m}}^{(j+1)2^{n-m}-1} \!\!\!\! C 2^{-n} |u(k2^{-n})|_H  +  \sum\limits_{\ell=n}^\infty \sum\limits_{r=j2^{\ell-m}}^{(j+1)2^{\ell-m}-1} \!\! |\varphi_{\ell+1,2r+1}\left(b_{n,\lfloor r2^{n-\ell}\rfloor} ; u\left((2r+1)2^{-\ell-1}\right), u\left(r2^{-\ell}\right)\right)|_H   \]

\[ = \!\!\! \sum\limits_{k=j2^{n-m}}^{(j+1)2^{n-m}-1} \!\!\!\!\!\! C 2^{-n} |u(k2^{-n})|_H  +  \sum\limits_{\ell=n}^\infty \sum\limits_{r=j2^{\ell+1-m}}^{(j+1)2^{\ell+1-m}-2} \!\!\!\!\!\! |\varphi_{\ell+1,r+1} \left(b_{n,\lfloor r2^{n-\ell-1}\rfloor} ; u\left((r+1) 2^{-\ell-1}\right), u\left(r 2^{-\ell-1}\right)\right)|_H . \]

We set for $\ell \geq n$

\[ \Lambda_\ell := \sum\limits_{r=j2^{\ell+1-m}}^{(j+1)2^{\ell+1-m}-2} \!\!   |\varphi_{\ell+1,r+1}\left(b_{n,\lfloor r2^{n-\ell-1}\rfloor} ; u\left((r+1)2^{-(\ell+1)}\right), u\left(r2^{-(\ell+1)}\right)\right)|_H   \]

and obtain

\[ \sum\limits_{k=j2^{n-m}}^{(j+1)2^{n-m}-1} \!\!\!\!\!\! \left| u((k+1)2^{-n}) - u(k2^{-n}) - \varphi_{n,k}(b_{n,k} ; u(k2^{-n})) \right|_H  \leq  \!\!\!\!\! \sum\limits_{r=j2^{\ell-m}}^{(j+1)2^{\ell-m}-1} \!\!\!\!\!\! C 2^{-n} |u(k2^{-n})|_H  +  \sum\limits_{\ell=n}^\infty \Lambda_\ell .  \label{T27E20} \]

From the reversed triangle inequality we deduce

\customlabel{THM-FINAL-MAIN}
\begin{align*}\tag{\ref{THM-FINAL-MAIN}}
\begin{split} \hspace{-80mm} &\sum\limits_{k=j2^{n-m}}^{(j+1)2^{n-m}-1} \!\!\!\!\!\!\!  |u((k+1)2^{-n}) - u(k2^{-n})|_H  \\  \hspace{20mm} \leq  &\sum\limits_{k=j2^{n-m}}^{(j+1)2^{n-m}-1}  \!\!\!\!\! \!\!\! \left( C 2^{-n}  |u(k2^{-n})|_H  +  |\varphi_{n,k}(b_{n,k} ; u(k2^{-n}))|_H \right)  +  \sum\limits_{\ell=n}^\infty \Lambda_\ell .
\end{split}
\end{align*}

The idea of the proof is the following: We will obtain estimates for the two sums on the right-hand side of the above inequality \eqref{THM-FINAL-MAIN}. For the first sum we simply use Theorem \ref{THM-SIGMA} (in the form of Corollary \ref{COR-SIGMA-RHO}) to obtain estimate \eqref{THM-FINAL-SIGMA}. We will split the second sum in the cases $\ell < N$ and $N \leq \ell$. In the first case we use Corollary \ref{COR-GLUING}, which will lead us to inequality \eqref{THM-FINAL-OMEGA-EULER}. For the second case we have to do a more direct computation, which heavily relies on the fact that $u$ is Lipschitz continuous (inequality \eqref{THM-FINAL-OMEGA-L}).\\

Combining all of this will result the final bound \eqref{THM-FINAL-OMEGA-FULL-ESTIMATE}. Using the knowledge of the already established estimate \eqref{THM-FINAL-MAIN} and the definition of $\alpha$ \eqref{THM-FINAL-DE-ALPHA} we will be able to estimate $\alpha$ in terms of $\beta$ (inequality \ref{THM-FINAL-ALPHA-ESTIMATE}). Feeding this back into inequality \eqref{THM-FINAL-DE-ALPHA} for $n=m$ completes the proof.\\

We will now estimate the two sums on the right-hand side starting with the $\varphi_{n,k}$ sum. We apply Corollary \ref{COR-SIGMA-RHO} to obtain

\[  \hspace{-19mm}  \sum\limits_{k=j2^{n-m}}^{(j+1)2^{n-m}-1}  \! \left( C 2^{-n} |u(k2^{-n})|_H  +  |\varphi_{n,k}(b_{n,k} ; u(k2^{-n}))|_H \right)  \]

\[ \leq  \sum\limits_{k=j2^{n-m}}^{(j+1)2^{n-m}-1} \!\!\! \left(C 2^{-n} |u(k2^{-n})|_H   +   C_\eps n^{\frac12+\frac1\gamma} 2^{-n/2} \left( \left| u(k2^{-n})\right|_H + 2^{-2^{n}} \right)\right) \]

\[ \hspace{-37mm}  \leq \sum\limits_{k=j2^{n-m}}^{(j+1)2^{n-m}-1} \!\!\! 2 \tilde C_\eps n^{\frac12+\frac1\gamma} 2^{-n/2} \left( \left| u(k2^{-n})\right|_H + 2^{-2^{n}} \right) \]

and since $n\geq m$ this is smaller than

\[ \hspace{-9mm}  2\tilde C_\eps n^{\frac12+\frac1\gamma} 2^{-n/2} \sum\limits_{k=j2^{n-m}}^{(j+1)2^{n-m}-1}  \left( |u(k2^{-n})|_H + 2^{-2^{m}} \right) \]

\[  = 2\tilde C_\eps n^{\frac12+\frac1\gamma} 2^{-n/2} \left( 2^{n-m} 2^{-2^{m}} +  \sum\limits_{k=j2^{n-m}}^{(j+1)2^{n-m}-1} |u(k2^{-n})|_H \right) .  \]

Again, using that $n\in\{m,...,N\}$ and the definition of $\psi_n$ (equation \eqref{THM-FINAL-DE-PSI}) this can be written as

\[ 2\tilde C_\eps n^{\frac12+\frac1\gamma} 2^{-n/2} \left( 2^{n-m} 2^{-2^{m}} +   \psi_n \right) . \]

Using inequality \eqref{THM-FINAL-PSI-ESTIMATE} this can be further estimated by

\[ 2\tilde C_\eps n^{\frac12+\frac1\gamma} 2^{-n/2} \left( 2^{n-m} 2^{-2^{m}}  +   2^{n-m} \left(\beta + \alpha 2^{-m/3} \right)\right) \]

\[ = 2\tilde C_\eps n^{\frac12+\frac1\gamma} 2^{n/2} 2^{-m} \left( 2^{-2^{m}}  +  \beta + \alpha 2^{-m/3} \right) \overset{2^{-2^{m}}\leq\beta}\leq  4\tilde C_\eps n^{\frac12+\frac1\gamma} 2^{n/2} 2^{-m} \left( \beta  +  \alpha 2^{-m/3} \right) \]

and hence for the first sum we obtain for all $n\in\{m,...,N\}$

\customlabel{THM-FINAL-SIGMA}
\[  \!\!\!\! \sum\limits_{k=j2^{n-m}}^{(j+1)2^{n-m}-1} \!\!\!\!\!\!\!\!  C 2^{-n} |u(k2^{-n})|_H  +  |\varphi_{n,k}(b_{n,k} ; u(k2^{-n}))|_H  \leq 4\tilde C_\eps n^{\frac12+\frac1\gamma} 2^{n/2} 2^{-m} \left( \beta  +  \alpha 2^{-m/3} \right)  .    \tag{\ref{THM-FINAL-SIGMA}}  \]

\newpage

Now consider the term $\Lambda_\ell$ for $\ell \geq N$. Applying Corollary \ref{COR-SIGMA-RHO} we obtain

\[ \hspace{-6mm}  \sum\limits_{\ell=N}^\infty \Lambda_\ell = \sum\limits_{\ell=N}^\infty  \sum\limits_{r=j2^{\ell+1-m}}^{(j+1)2^{\ell+1-m}-2}  \left|\varphi_{\ell+1,r+1}\left(b_{n,\lfloor r2^{n-\ell-1}\rfloor} ; u\left((r+1)2^{-(\ell+1)}\right), u\left(r2^{-(\ell+1)}\right)\right)\right|_H  .  \]

\[ \hspace{15mm}  \leq \sum\limits_{\ell=N}^\infty  \sum\limits_{r=j2^{\ell+1-m}}^{(j+1)2^{\ell+1-m}-2}  C_\eps \left( \sqrt{\ell+1} 2^{-\ell/6} \left| u\left((r+1)2^{-(\ell+1)} \right) - u\left(r2^{-(\ell+1)}\right) \right|_\infty  +  2^{-2\ell} \right) .  \]

By the Lipschitz continuity of $u$ this is smaller than

\[ C_\eps \sum\limits_{\ell=N}^\infty  \sqrt{\ell+1} 2^{-\ell/6} \sum\limits_{r=j2^{\ell+1-m}}^{(j+1)2^{\ell+1-m}-1} \left( |(r+1)2^{-\ell-1} - r2^{-\ell-1}|  +  2^{-\ell} \right) \]

\[ \hspace{-37mm}  =  C_\eps \sum\limits_{\ell=N}^\infty \sqrt{\ell+1} 2^{-\ell/6} \sum\limits_{r=j2^{\ell+1-m}}^{(j+1)2^{\ell+1-m}-1}  ( 2^{-\ell} + 2^{-\ell-1} ) \]

\[ = \frac32 C_\eps \sum\limits_{\ell=N}^\infty 2^{\ell-m} \sqrt{\ell+1} 2^{-\ell/6} 2^{-\ell} = \frac32 C_\eps 2^{-m} \sum\limits_{\ell=N}^\infty  \sqrt{\ell+1} 2^{-\ell/6} \leq 2 C_\eps 2^{-m} 2^{-N/7} . \]

And hence we obtain

\customlabel{THM-FINAL-OMEGA-L}
\[  \sum\limits_{\ell=N}^\infty \Lambda_\ell \leq 2C_\eps 2^{-m} 2^{-N/7} .   \tag{\ref{THM-FINAL-OMEGA-L}}  \]

Now consider the case $n \leq \ell \leq N$. We define

\[  \gamma_{\ell,r} := u((r+1)2^{-\ell}) - u(r2^{-\ell}) - \varphi_{\ell,r}(b_{n,\lfloor r2^{n-\ell-1}\rfloor} ; u(r2^{-\ell})), \qquad\qquad \forall r\in\{0,...,2^\ell-1\}  \]

and note that due to inequality \eqref{THM-FINAL-MAIN} we have

\customlabel{THM-FINAL-GAMMA-SMALLER-OMEGA}
\[ \sum\limits_{r=j2^{\ell-m}}^{(j+1)2^{\ell-m}-1}  \!\! |\gamma_{\ell,r}|_H  \leq \sum\limits_{\ell'=\ell}^\infty \Lambda_{\ell'} .   \tag{\ref{THM-FINAL-GAMMA-SMALLER-OMEGA}}    \]

Recall the definition of $\Lambda_\ell$:

\[ \Lambda_\ell = \sum\limits_{r=j2^{\ell+1-m}}^{(j+1)2^{\ell+1-m}-2} \!\!  \left|\varphi_{\ell+1,r+1}\left(b_{n,\lfloor r2^{n-\ell-1}\rfloor} ; u\left((r+1)2^{-(\ell+1)}\right), u\left(r2^{-(\ell+1)}\right)\right)\right|_H  . \]

Using Corollary \ref{COR-GLUING} yields that this is bounded from above by

\[   C_\eps \left[ 2^{-\ell} \psi_\ell + 2^{-\ell/24} \!\! \sum\limits_{r=j2^{\ell-m}}^{(j+1)2^{\ell-m}-1} \!\!\!\! |\gamma_{\ell,r}|_H  +  2^{-3\ell/4} |u(j2^{-m})|_H  +  2^{\ell+1-m} 2^{-2^{\theta \ell}} \right] .  \]

Since $\ell \leq N$ we can use inequality \eqref{THM-FINAL-PSI-ESTIMATE} and the assumption $|u(j2^{-m})|_H \leq \beta$ to obtain the following estimate

\[ \sum\limits_{\ell=n}^N \Lambda_\ell  \leq  C_\eps \sum\limits_{\ell=n}^N \left[ 2^{-m} (\beta + \alpha 2^{-m/3})  +  2^{-\ell/24} \!\! \sum\limits_{r=j2^{\ell+1-m}}^{(j+1)2^{\ell+1-m}-1} \!\!\! |\gamma_{\ell,r}|_H  + 2^{-3\ell/4} \beta +  2^{\ell+1-m} 2^{-2^{\theta \ell}} \right] .  \]

Using inequality \eqref{THM-FINAL-GAMMA-SMALLER-OMEGA} and $2^\ell 2^{-2^{\theta\ell}} \leq 2^{n-\ell} 2^n 2^{-2^{\theta n}}$ this can be further estimated by

\[ \hspace{4mm}  C_\eps \left[ 2^{-m} N (\beta + \alpha 2^{-m/3})  +  \sum\limits_{\ell=n}^N 2^{-\ell/24} \sum\limits_{\ell'=\ell}^\infty \Lambda_{\ell'}  + 2^{-3n/4} \beta +  2^{-m} 2^{n+1} 2^{-2^{\theta n}} \right]  \]

and since $m \leq n$ this is smaller than

\[  36 C_\eps \left[ 2^{-m} N (\beta + \alpha 2^{-m/3})  +  2^{-m/24} \sum\limits_{\ell=n}^\infty \Lambda_\ell  + 2^{-3m/4} \beta +  2 \cdot 2^{-2^{\theta m}} \right] . \]

Recall that by \eqref{THM-FINAL-DE-N} we have $2^{-3m/4} = 2^{-m} 2^{m/4} \leq 2^{-m}N$ as well as $2^{-2^{\theta m}} \leq 2^{-m} \beta$, so that in conclusion we deduce

\customlabel{THM-FINAL-OMEGA-EULER}
\[ \sum\limits_{\ell=n}^N \Lambda_\ell \leq 144 C_\eps 2^{-m} N (\beta  + \alpha 2^{-m/3}) + \frac12 \sum\limits_{\ell=n}^\infty \Lambda_\ell ,  \tag{\ref{THM-FINAL-OMEGA-EULER}}  \]

where we have used that $36 C_\eps 2^{-m/24} \leq 36 C_\eps  2^{-m_0/24} \leq \frac12$.

\medskip

Putting together the both estimates \eqref{THM-FINAL-OMEGA-EULER} and \eqref{THM-FINAL-OMEGA-L} for $\Lambda_\ell$ we have

\[ \sum\limits_{\ell=n}^\infty \Lambda_\ell = \sum\limits_{\ell=n}^N \Lambda_\ell + \sum\limits_{\ell=N+1}^\infty \Lambda_\ell  \leq 144 C_\eps 2^{-m} N (\beta  + \alpha 2^{-m/3}) + \frac12 \sum\limits_{\ell=n}^\infty \Lambda_\ell  +   2 C_\eps 2^{-m} 2^{-N/7} . \]

Henceforth, we deduce

\[ \sum\limits_{\ell=n}^\infty \Lambda_\ell \leq 288 C_\eps 2^{-m} N (\beta  + \alpha 2^{-m/3}) +   4 C_\eps 2^{-m} 2^{-N/7}  \]

and since $N$ reads $N = 7 \lfloor \log_2(1/\beta) \rfloor$ this expression is bounded by

\[ 292 C_\eps 2^{-m} N (\beta  + \alpha 2^{-m/3}) .  \]

Therefore, we have

\customlabel{THM-FINAL-OMEGA-FULL-ESTIMATE}
\[ \sum\limits_{\ell=n}^\infty \Lambda_\ell  \leq  292 C_\eps 2^{-m} N (\beta  + \alpha 2^{-m/3})  .  \tag{\ref{THM-FINAL-OMEGA-FULL-ESTIMATE}}  \]

Looking back to inequality \eqref{THM-FINAL-MAIN}, with the help of \eqref{THM-FINAL-SIGMA} and \eqref{THM-FINAL-OMEGA-FULL-ESTIMATE}, we estimate the sum by

\[ \!\! \sum\limits_{k=j2^{n-m}}^{(j+1)2^{n-m}-1} \!\!\!\!\!\!\! |u((k+1)2^{-n}) - u(k2^{-n})|_H \leq  \!\!\!\!\!\!\! \sum\limits_{k=j2^{n-m}}^{(j+1)2^{n-m}-1} \!\!\!\!\!\!\!\! \left( C2^{-n} |u(k2^{-n})|_H  +  |\varphi_{n,k}(b_{n,k} ; u(k2^{-n}))|_H \right)  +  \sum\limits_{\ell=n}^\infty \Lambda_\ell  \]

\[ \leq 4 \tilde C_\eps  n^{\frac12+\frac1\gamma} 2^{n/2} 2^{-m} \left[ \beta  +  \alpha 2^{-m/3} \right]  +  292 C_\eps 2^{-m} N \left[ \beta  + \alpha 2^{-m/3} \right ] \]

\[ \hspace{-30mm}  \leq 292 C_\eps' 2^{-m} \left[ n^{\frac12+\frac1\gamma} 2^{n/2}  +  N \right] \cdot \left[ \beta  +  \alpha 2^{-m/3} \right] . \]

Note that the above argument holds for all $n \in\{ m , ..., N\}$. Hence, by the minimality of $\alpha$ and inequality \eqref{THM-FINAL-DE-ALPHA} we have

\[ \alpha 2^{-m}  \left[ n^{\frac12+\frac1\gamma} 2^{n/2}  +  N \right] \leq 292 C_\eps' 2^{-m}  \left[ n^{\frac12+\frac1\gamma} 2^{n/2}  +  N \right] \cdot \left[ \beta  +  \alpha 2^{-m/3} \right] \]

for all $n\in\{m,...,N\}$. This implies that

\[ \alpha \leq 292 C_\eps' \left[ \beta  +  \alpha 2^{-m/3} \right] .  \]

Since $292 C_\eps' 2^{-m/3} \leq 292 C_\eps' 2^{-m_0/3} \leq \frac12$ holds for all $m\geq m_0$. It now follows

\[ \alpha \leq 292 C_\eps' \beta  +  \alpha 292 C_\eps' 2^{-m/3} \leq 292 C_\eps' \beta + \frac\alpha2 . \]

From which we deduce that

\customlabel{THM-FINAL-ALPHA-ESTIMATE}
\[ \alpha \leq 584 C_\eps' \beta .  \tag{\ref{THM-FINAL-ALPHA-ESTIMATE}}  \]

Setting $n=m$ in \eqref{THM-FINAL-DE-ALPHA} reads

\[  |u((j+1)2^{-m}) - u(j2^{-m})|_H \!\! \overset{\eqref{THM-FINAL-DE-ALPHA}}\leq \!\! \alpha 2^{-m} \left[ m^{\frac12+\frac1\gamma} 2^{m/2} + N \right] .  \]

Putting $|u(j2^{-m})|_H$ to the right-hand side yields

\[ |u((j+1)2^{-m})|_H   \leq |u(j2^{-m})|_H + \alpha 2^{-m} \left[ m^{\frac12+\frac1\gamma}  2^{m/2} + N \right]   \]

and since we have $\alpha \leq 584 C_\eps' \beta$ as well as $m^{\frac12+\frac1\gamma} 2^{m/2} \leq N$ by definition of $N$ using our estimate \eqref{THM-FINAL-ALPHA-ESTIMATE} yields that the above expression is smaller than

\[  \beta + 584 C_\eps' \beta 2^{-m} N  = \beta \left( 1 + 584 C_\eps' 2^{-m} \lfloor \log_2(1/\beta) \rfloor \right)  \leq \beta \left( 1 + K 2^{-m} \log_2(1/\beta) \right) ,  \]

where the constant is defined as $K := 584 C_\eps'$ which completes the proof.

\qed

\end{prof}

\begin{cor}
\xlabel{COR-FINAL}

Let $f$ be a $H$-valued Borel function such that the Assumption \ref{ASS} is fulfilled then there exists a set $N_f\subseteq\Omega$ with $\P[N_f]=0$ such that for all $\omega\in N_f^c$ if $u$ is a solution to

\[ u(t) = \int\limits_0^t e^{-(t-s)A} ( f(s, u(s) + Z^A_s(\omega)) - f(s, Z^A_s(\omega)) ) \d s , \qquad    \forall t\in [0,1] . \]

then $u \equiv 0$.

\end{cor}

\begin{prof}

\textbf{Step 1:}\\
\\
Let $0 < \eps < \frac1{40}$ and $\Omega_{\eps,f}$ be the of set of Theorem \ref{THM-FINAL}.
Fix $\omega \in \Omega_{\eps,b}$ and let $u$, as stated in the assertion, be a solution to the above equation. Since $\|f\|_\infty \leq 1$ the function $u$ is Lipschitz continuous with Lipschitz constant at most $2$.
Furthermore, Assumption \ref{ASS} on $f$ implies that $u$ is $Q$ as well as $Q^A$-valued.

\medskip

Therefore $u\in\Phi$. Applying Theorem \ref{THM-FINAL} gives us a $K>0$ and $m_0\in\N$. For sufficiently large $m\in\N$ (i.e.~$K \leq \ln(2) 2^m$ and $m\geq m_0$) we define

\[ \beta_0  := 2^{m-2^{\theta m}} \]

and

\[ \beta_{j+1} := \beta_j(1+K2^{-m} \log_2(1/\beta_j)) \]

for $j \in \{0, ..., 2^m -1 \}$. By the very definition we have

\[ \beta_{j+1} - \beta_j = K2^{-m} \beta_j \log_2(1/\beta_j)  \]

for every $j \in \{0, ..., 2^m -1 \}$. Hence, Lemma \ref{LEM-GRONWALL} is applicable which implies that

\[ \beta_j \leq \exp\left( \log_2(\beta_0) e^{-2K-1} \right) = \exp\left( \left(m-2^{\theta m}\right) e^{-2K-1} \right) \]

\[ \hspace{-17mm}  \leq \exp\left( - \ln(2) 2^{\left(\frac\theta2+\frac14\right) m} \right) = 2^{-2^{\left(\frac\theta2+\frac14\right) m}} . \]

Together with the fact that $\beta_j$ is increasing we have

\[ 2^{m-2^{\theta m}} \leq \beta_j \leq 2^{-2^{(\frac\theta2+\frac14) m}} .  \]

Since $u$ is a solution to equation \eqref{DE-U} we know that $u(0) = 0 \leq \beta_0$, so that we are able to invoke Theorem \ref{THM-FINAL} from which we deduce that $|u(2^{-m})|_H \leq \beta_1$. By induction on $j$ we obtain

\[ |u(j2^{-m})|_H \leq \beta_j \leq \beta_{2^m} \leq 2^{-2^{\left(\frac\theta2+\frac14\right) m}} , \qquad \forall j\in\{0, ..., 2^m \} \]

for every $j \in \{0, ..., 2^m -1 \}$. By letting $m\rightarrow\infty$, we deduce that $u$ vanishes at all dyadic points.
By continuity of $u$ it follows $u \equiv 0$.\\

\textbf{Step 2:}\\
\\
Let $k\in\N$. By setting $\eps:=1/k$ in Step $1$ we conclude that there is $\Omega_{k,f} \subseteq\Omega$ with $\P[\Omega_{k,f}^c] \leq 1/k$ such that $u \equiv 0$ for all $\omega\in \Omega_{k,b}$. By defining

\[ N_f := \bigcap\limits_{k=41}^\infty \Omega_{k,f} \]

we have $u\equiv 0$ for all $\omega\in N_f^c$ which concludes the proof.

\qed

\end{prof}

\begin{prof}[of the main result \ref{PRO-MAIN}]

By Corollary \ref{COR-FINAL} the assumption of Proposition \ref{PRO-GIRSANOV} is fulfilled, so that by invoking Proposition \ref{PRO-GIRSANOV} the conclusion of Proposition \ref{PRO-MAIN} follows. Theorem \ref{THM-MAIN} then follows from Proposition \ref{PRO-MAIN} as explained in the introduction.

\qed

\end{prof}

\newpage


\begin{footnotesize}

\end{footnotesize}


\begin{thebibliography}{BFGM14]}

\bibitem[BM16]{BM16}
O.~Butkovsky, L.~Mytnik.
\newblock {\textit{Regularization by noise and flows of solutions for a stochastic heat equation}}.
\newblock {ArXiv e-prints}, 2016.
\newblock \href{http://arxiv.org/abs/1610.02553v2}{arXiv: 1610.02553v2}.

\bibitem[Bog98]{Bog98}
V.~I.~Bogachev.
\newblock {\textit{Gaussian Measures}}.
\newblock Mathematical Surveys and Monographs. American Math. Soc., 1998.

\bibitem[BFGM14]{BFGM14}
L.~Beck, F.~Flandoli, M.~Gubinelli, M.~Maurelli.
\newblock {\textit{Stochastic ODEs and stochastic linear PDEs with critical
drift: regularity, duality and uniqueness.}}
\newblock {ArXiv e-prints}, 2014.
\newblock \href{http://arxiv.org/abs/1401.1530}{arXiv: 1401.1530}.

\bibitem[Dav07]{Dav07}
A.~M.~{Davie}.
\newblock {\textit{Uniqueness of solutions of stochastic differential equations}}.
\newblock {ArXiv e-prints}, 2007.
\newblock \href{http://arxiv.org/abs/0709.4147}{arXiv: 0709.4147}.

\bibitem[Dav76]{Dav76}
B.~Davis.
\newblock {\textit{On the $l^p$ norms of stochastic integrals and other martingales}}.
\newblock {Duke Math.~J.}, \textbf{43} no.~4, pp.~697--704, 1976.

\bibitem[DFPR13]{DFPR13}
G.~Da Prato, F.~Flandoli, E.~Priola, M.~R{\"o}ckner.
\newblock {\textit{Strong uniqueness for stochastic evolution equations in Hilbert spaces perturbed by a bounded measurable drift}}.
\newblock {Ann.~Probab.} \textbf{41}, no.~5, pp.~3306--3344, 2013.

\bibitem[Dug51]{Dug51}
J.~Dugundji.
\newblock {\textit{An extension of Tietze's theorem}}.
\newblock {Pac.~J.~Math.}, \textbf{1}, no.~3, pp.~353--367, 1951.

\bibitem[DZ92]{DZ92}
G.~Da Prato, J.~Zabczyk.
\newblock{\textit{Stochastic Equations in Infinite Dimensions}}.
\newblock{Cambridge Univ.~Pr.},
pp.~454, 1992.

\bibitem[Fla10]{Fla10}
F.~Flandoli.
\newblock{\textit{Random Perturbation of PDEs and Fluid Dynamic Models}}.
\newblock{\'Ecole d'\'Et\'e de Probabilit\'es de Saint-Flour XL},
pp.~182, 2010.

\bibitem[Kal97]{Kal97}
O.~Kallenberg.
\newblock{\textit{Foundations of Modern Probability}}.
\newblock{Springer New York. Probability and its Applications: A Series of the Applied Probability Trust},
1997.

\bibitem[Kho14]{Kho14}
D.~Khoshnevisan.
\newblock{\em Analysis of stochastic partial differential equations}
\newblock{CBMS regional conference series in mathematics} \textbf{119}.
\newblock Am.~Math.~Soc., 2014.

\bibitem[Ose12]{Ose12}
A.~Osekowski.
\newblock {\em Sharp Martingale and Semimartingale Inequalities},
{Monografie Matematyczne} \textbf{72}
\newblock Springer, 2012.

\bibitem[Pri15]{Pri15}
E.~Priola.
\newblock {\textit{Davie's type uniqueness for a class of SDEs with jumps}}.
\newblock { ArXiv e-prints}, 2015.
\newblock \href{http://arxiv.org/abs/1509.07448v2}{arXiv: 1509.07448v2}.

\bibitem[LR15]{LR15}
W.~Liu, M.~R{\"o}ckner.
\newblock {\em Stochastic Partial Differential Equations: An Introduction}.
\newblock Springer International Publishing, pp.~266, 2015.

\bibitem[RSZ08]{RSZ08}
M.~R{\"o}ckner, B.~Schmuland, X.~Zhang.
\newblock {\em Yamada--{W}atanabe {T}heorem for stochastic evolution equations in infinite dimensions}.
\newblock {Cond.~Matt.~Phys.}, \textbf{11} no.~2, pp.~247--259, 2008.

\bibitem[Sha14]{Sha14}
A.~{Shaposhnikov}.
\newblock {\textit{Some remarks on Davie's uniqueness theorem}}.
\newblock {ArXiv e-prints}, 2014.
\newblock \href{http://arxiv.org/abs/1401.5455}{arXiv: 1401.5455}.

\bibitem[Wre16]{Wre16}
L.~Wresch.
\newblock {\textit{An exponential estimate for Hilbert space-valued Ornstein--Uhlenbeck processes}}.
\newblock {ArXiv e-prints}, 2016.
\newblock \href{http://arxiv.org/abs/1612.07745}{arXiv: 1612.07745}.

\end{thebibliography}
\end{document}